%% file: IntAddr2012.tex
\documentclass[12pt,a4paper,reqno]{amsart}

\usepackage{graphicx}
\usepackage{latexsym}
\usepackage{amsfonts}
\usepackage{epsfig}

\makeindex

\include{TreesInclude}

\newenvironment{changemargin}[2]{%
  \begin{list}{}{%
    \setlength{\topsep}{0pt}%
    \setlength{\leftmargin}{#1}%
    \setlength{\rightmargin}{#2}%
    \setlength{\listparindent}{\parindent}%
    \setlength{\itemindent}{\parindent}%
    \setlength{\parsep}{\parskip}%
  }%
  \item[]}{\end{list}}

\def\reminder#1{{\sf #1} }
\def\hide#1{}

\renewcommand{\subsection}{\section}

\newcommand{\version}[2] {#2} 

\newcommand{\comment}[1]{\marginpar{#1}}

\renewcommand{\comment}[1]{}

\newcommand{\CorNumberEmbed}{Corollary~{\bf 4.11}}

\newcommand{\PropExactPeriod}{Proposition~{\bf 5.16}}
\newcommand{\CorMaxShiftAdmiss}{Corollary~{\bf 5.20}}
\newcommand{\PropIntAddr}{Proposition~{\bf 6.8}}

\newcommand{\SecRenorm}{Section~{\bf 10}}
\newcommand{\ThmCorrespondencence}{Theorem~{\bf 9.4}}
\newcommand{\ThmNoGhostLimbs}{Theorem~{\bf 9.12}}
\newcommand{\CorIrrationalIntAngles}{Corollary~{\bf 9.13}}
\newcommand{\CorBoundaryRenormalization}{Corollary~{\bf 9.27}} 
\newcommand{\PropWidthWake}{Proposition~{\bf 9.34}} 
\newcommand{\LemPeriodicItineraryUpper}{Lemma~{\bf 14.6}}
\newcommand{\LemExactPeriod}{Lemma~{\bf 19.2}}
\newcommand{\SecAdmissTwice}{Sections~{\bf 5} and {\bf 14}}
\newcommand{\LemNumberEmbedTwo}{Lemma~{\bf 16.9}}

\newcommand{\dit}{\texttt {d}}

\begin{document}
\title[Internal Addresses and Galois Groups
]
{Internal Addresses of the Mandelbrot Set \\and Galois Groups of Polynomials
}
\author{Dierk Schleicher}

\begin{abstract}
We describe an interesting interplay between symbolic dynamics, the structure of the Mandelbrot set, permutations of periodic points achieved by analytic continuation, and Galois groups of certain polynomials.

Internal addresses are a convenient and efficient way of describing the combinatorial structure of the Mandelbrot set, and of giving geometric meaning to the ubiquitous kneading sequences in human-readable form (Sections~\ref{Sub:IntAddr} and \ref{Sub:GeometryIntAddr}). A simple extension, \emph{angled internal addresses}, distinguishes combinatorial classes of the Mandelbrot set and in particular distinguishes hyperbolic components in a concise and dynamically meaningful way.

This combinatorial description of the Mandelbrot set makes it possible to derive existence theorems for certain kneading sequences and internal addresses in the Mandelbrot set (Section~\ref{Sec:Permutations}) and to give an explicit description of the associated parameters. These in turn help to establish some algebraic results about permutations of periodic points and to determine Galois groups of certain polynomials (Section~\ref{Sub:Permutations}). 

Through internal addresses, various areas of mathematics are thus related in this manuscript, including symbolic dynamics and permutations, combinatorics of the Mandelbrot set, and Galois groups.
\end{abstract}

\maketitle

\renewcommand{\thefootnote}{}
\footnotetext{\emph{Address:} Jacobs University, Research I, Postfach 750 561, D-28725 Bremen, Germany; dierk$@$jacobs-university.de}
\renewcommand{\thefootnote}{\arabic{footnote}}

\begin{changemargin}{20mm}{20mm}
\tableofcontents
\end{changemargin}

\input{SecIntAddr}

\input{RefsIntAddr}



\end{document}

%% file: TreesInclude.tex
\newcounter{null}
\def\Nzero{{\mathbb N}}
\def\Q{\mathbb {Q}}
\def\Z{\mathbb {Z}}
\def\R{\mathbb {R}}
\def\C{\mathbb {C}}

\newcommand{\Circle}{{\mathbb S}^1}

\newcommand{\heading}[1]{\par\medskip {\bf #1.}} 

\newfont\bbf{msbm10 at 12pt}
\newcommand\ovl[1]{\overline {#1}}

\newfont\script{eusm10 at 12pt}
\newcommand\M{\text{\script M}}

\newcommand{\wake}[1]{\mathcal{W}_{{#1}}}

\newtheorem{satz}{\hspace{-1.5mm}}[section]

\def\lineclear
      {\rule{0pt}{0pt}\nopagebreak\par\nopagebreak\noindent}
\newenvironment {lemma} [1]
{\begin{satz} {\bf Lemma (#1) }\nopagebreak}{\end{satz}}
\newenvironment {defi} [1]
{\begin{satz} {\bf Definition (#1)}\nopagebreak}{\end{satz}}
\newenvironment {theo} [1]
{\begin{satz} {\bf Theorem (#1)}\nopagebreak}{\end{satz}}
\newenvironment {prop} [1]
{\begin{satz} {\bf Proposition (#1) }\nopagebreak}{\end{satz}}
\newenvironment {coro} [1]
{\begin{satz} {\bf Corollary (#1) }\nopagebreak}{\end{satz}}
\newenvironment {coro*}
{\begin{satz} {\bf Corollary }\nopagebreak}{\end{satz}}

\newenvironment {algo} [1]
{\begin{satz} {\bf Algorithm (#1) }\nopagebreak}{\end{satz}}

\newenvironment {conjecture} [1]
{\begin{satz} {\bf Conjecture (#1) }\nopagebreak}{\end{satz}}

\def\proof{\medskip\par\noindent{\sc Proof. }}
\def\proofof#1{\medskip\par\noindent{\sc Proof of #1}.}

\newcommand\nix{\rule{0pt}{2pt}}
\renewcommand{\qed}{\qedd\par\medskip\noindent}
\newcommand{\qedd}{\nix\nolinebreak\hfill\hfill\nolinebreak$\Box$}

\def\remark{\par\medskip \noindent {\sc Remark. }}
\newcommand\Intro[1] {{\small #1}\par\medskip}

\def\reminder#1{{\sf #1} }
\def\hide#1{}

\def\0{{\tt 0}}
\def\1{{\tt 1}}
\def\*{{\tt \star}}
\newcommand\A{{\mathcal A}}
\newcommand\Abar{\ovl{\A}}
\def\IntAddr{\to}
\def\IntAdr{\to}

\renewcommand{\text}[1]{\mbox{#1}}  
\newcommand{\ph}{\varphi}
\renewcommand\phi{\varphi}

\renewcommand\theta{\vartheta}
\newcommand{\eps}{\varepsilon}

\newcommand{\orb}{\mbox{\rm orb}}

\newcommand{\disk}{\mbox{\bbf D}}
\newcommand{\diskbar}{\ovl\disk}
\newcommand{\RP}{P} 

\newcommand\sm{\setminus}

\newcounter{stepcount}
    {\begin{list}{Step~\arabic{stepcount}.}{\usecounter{stepcount}}}%
    {\end{list}}



%% file: SecIntAddr.tex


\Intro{
\version{
In the early 1990's, Devaney asked the question: {\em how can you tell where in the Mandelbrot set a given external ray lands, without having Adrien Douady at your side?} 
In this section, we provide an answer to this question in terms of internal addresses: these are a convenient and efficient way of describing the combinatorial structure of the Mandelbrot set, and of giving geometric meaning to the ubiquitous kneading sequences in human-readable form. 

We give the most important results
from \cite{IntAdr} about the combinatorial geometry of internal addresses. These will be used in Section~\ref{Sec:Permutations} to derive some algebraic results about permutations of periodic points and Galois groups.
}
{
\hide{
We describe an interesting interplay between symbolic dynamics, the structure of the Mandelbrot set, permutations of periodic points achieved by analytic continuation, and Galois groups of certain polynomials.
In the early 1990's, Devaney asked the question: {\em how can you tell where in the Mandelbrot set a given external ray lands, without having Adrien Douady at your side?} 
We provide an answer to this question in terms of \emph{internal addresses}: these are a convenient and efficient way of describing the combinatorial structure of the Mandelbrot set, and of giving geometric meaning to the ubiquitous kneading sequences in human-readable form. A simple extension, \emph{angled internal addresses}, distinguishes combinatorial classes of the Mandelbrot set and in particular distinguishes hyperbolic components in a concise and dynamically meaningful way.
This combinatorial description of the Mandelbrot set makes it possible to derive quite easily existence theorems for certain kneading sequences and internal addresses in the Mandelbrot set; these in turn help to establish some algebraic results about permutations of periodic points and to determine Galois groups of certain polynomials. Many results in this paper were first announced in \cite{IntAdr}.
}
}
}

\section{Introduction}

The combinatorial structure of the Mandelbrot set has been studied by many people, notably in terms of \emph{quadratic minor laminations} by Thurston \cite{thurston}, \emph{pinched disks} by Douady \cite{DoCompacts}, or \emph{orbit portraits} by Milnor \cite{MiOrbits}\version{ (see also Appendix~\ref{appendix})}{}. All of these results are modeled on parameter rays of the Mandelbrot set at periodic angles, as well as their landing properties: parameter rays at periodic angles of fixed period are known to land together in pairs, so these pairs subdivide the complex plane into finitely many components. Internal addresses are a natural way to distinguish these components, and thus to describe the combinatorial structure of the Mandelbrot set in an efficient way.

One feature of internal addresses is that they provide a good language to describe of the combinatorial structure of $\M$ (Sections~\ref{Sub:IntAddr}, \ref{Sub:GeometryIntAddr}, and \ref{Sub:ProofsIntAddr}). Internal addresses also allow us to decide which combinatorial data are realized in $\M$: we describe a particularly relevant kind of hyperbolic components that we call \emph{purely narrow}, and we give a necessary and sufficient description of the associated combinatorics (Section~\ref{Sec:Permutations}). 

Finally, we address an apparently rather different topic: when periodic points are continued analytically across parameter space along closed loops, which permutations can be so achieved? This question can be interpreted as the determination of certain Galois groups of dynamically defined polynomials, or by the structure of the ramified cover over parameter space given by periodic points of given period. We explicitly determine these Galois groups in Section~\ref{Sub:Permutations} and show how this question leads very naturally to questions about the combinatorial structure of $\M$ of exactly the kind that was answered in previous sections. 

Many results in this paper were first announced in \cite{IntAdr}. The Galois groups had been determined earlier in a more algebraic way by Bousch \cite{Bousch} and more recently by Morton and Patel \cite{MP}.

\section{Background}

Much of the combinatorial structure of the Mandelbrot set $\M$ is described in terms of parameter rays; these are defined as follows\version{ (see Section~\ref{SecMandelbrot})}{}. By \cite{Orsay}, the Mandelbrot set is compact, connected and full, i.e., there is a conformal isomorphism $\Phi\colon(\C\sm\M)\to(\C\sm\diskbar)$; it can be normalized so that $\Phi(c)/c\to 1$ as $c\to\infty$. The \emph {parameter ray} of $\M$ at angle $\theta\in\R/\Z$ is then defined as $R(\theta):=\Phi^{-1}(e^{2\pi i\theta}(1,\infty))$. The parameter ray $R(\theta)$ is said to \emph{land} at some point $c\in\partial\M$ if $\lim_{r\searrow 1} \Phi^{-1}(re^{2\pi i\theta})$ exists and equals $c$.

If $R(\theta)$ and $R(\theta')$ land at a common point $z$, we call $R(\theta)\cup R(\theta')\cup\{z\}$ a \emph{parameter ray pair} and denote it $\RP(\theta,\theta')$. Then $\C\sm\RP(\theta,\theta')$ has two components, and we say that $\RP(\theta,\theta')$ \emph{separates} two points or sets if they are in different components.
If a polynomial $p_c(z):=z^2+c$ has connected Julia set, then there are analogous definitions of \emph{dynamic rays} $R_c(\theta)$ and \emph{dynamic ray pairs} $\RP_c(\theta,\theta')$. General background can be found in Milnor \cite{MiBook}.

A \emph{hyperbolic component} of $\M$ is a connected component of the set of $c\in\M$ for which the map $p_c\colon z\mapsto z^2+c$ has an attracting periodic point of some period, say $n$. This period is necessarily constant throughout the component, and we say that $n$ is the period of the hyperbolic component. Every hyperbolic component $W$ has a distinguished boundary point called the \emph{root} of $W$; this is the unique parameter on $\partial W$ on which the attracting orbit becomes indifferent with multiplier $1$.

The most important properties of parameter rays at periodic angles are collected in the following well known theorem; see \version{Theorems~\ref{Thm:ParaRayPairs} and \ref{Thm:Correspondence} and Lavaurs' Lemma~\ref{Lem:Lavaurs}}{\cite{Orsay}, \cite{MiOrbits}, \cite{Fibers2}}.

\begin{theo}{Properties of Periodic Parameter Ray Pairs}
\label{Thm:PropertiesParaRayPairs} \lineclear
For every $n\ge 1$, every parameter ray $R(\theta)$ of period $n$ lands at the root of a hyperbolic component $W$ of period $n$, and the root of every hyperbolic component of period $n$ is the landing point of exactly two parameter rays, both of period $n$. 

If $R(\theta')$ is the other parameter ray landing at the same root as $R(\theta)$, then the ray pair $\RP(\theta,\theta')$ partitions $\C$ into two open components; let $\wake{W}$ be the component containing $W$: this is the \emph{wake of $W$}, and it does not contain the origin. This wake is the locus of parameters $c\in\C$ for which the dynamic rays $R_c(\theta)$ and $R_c(\theta')$ land together at a repelling periodic point; the ray pair $\RP_c(\theta,\theta')$ is necessarily characteristic.
\end{theo}

Here a dynamic ray pair is \emph{characteristic} if it separates the critical value from the rays $R_c(2^k\theta)$ and $R_c(2^k\theta')$ for all $k\ge 1$ (except of course from those on the ray pair $\RP_c(\theta,\theta')$ itself).

For period $n=1$, there is a single parameter ray $R(0)=R(1)$; the statement of the theorem holds if we count $0$ and $1$ separately; we do the same for dynamic rays.

In particular, we have the following important result.
\begin{coro}{The Ray Correspondence Theorem} 
\label{Cor:RayCorrespondence} \lineclear
For every parameter $c\in\C$, there is an angle-preserving bijection between parameter ray pairs separating the parameter $c$ from the origin and characteristic dynamic ray pairs in the dynamical plane of $p_c$ landing at repelling periodic points and separating the critical value $c$ from the origin.
\end{coro}

Parameter ray pairs at periodic angles have the following important property.
\begin{lemma}{Lavaurs' Lemma \cite{Lavaurs}}
\label{Lem:Lavaurs} \lineclear
If any parameter ray pair at periodic angles separates another parameter ray pair at equal period from the origin, then these ray pairs are separated by a parameter ray pair of strictly lower period. 
\end{lemma}

The following combinatorial concept will play an important role in our discussion.

\begin{defi}{Kneading Sequence}
\label{Def:Kneading} \lineclear
Every angle $\theta\in\Circle:=\R/\Z$ has an associated \emph{kneading sequence} $\nu(\theta)=\nu_1\nu_2\nu_3\dots$\version{ (Definition~\ref{DefKneading})}{}, defined as the itinerary of $\theta$ (under angle doubling) on the unit circle $\Circle$ with respect to the partition $\Circle\sm\{\theta/2,(\theta+1)/2\}$, so that 
\[
\nu_k=\left\{ \begin{array}{ll}
\1 & \text{ if }2^{k-1}\theta \in(\theta/2,(\theta+1)/2); \\
\0 & \text{ if }2^{k-1}\theta \in ((\theta+1)/2,\theta/2)=\Circle\sm  [\theta/2,(\theta+1)/2]; \\
\* & \text{ if }2^{k-1}\theta \in\{\theta/2,(\theta+1)/2\}.
\end{array}\right.
\]
\end{defi}
\hide{
$\nu_k=\1$ if $2^{k-1}\theta\in(\theta/2,(\theta+1)/2)$, $\nu_k=\0$ if $2^{k-1}\theta\in((\theta+1)/2,\theta/2)=\Circle\sm[\theta/2,(\theta+1)/2]$ and  $\nu_k=\*$ if $2^{k-1}\theta\in\{\theta/2,(\theta+1)/2\}$.
} 
If $\theta$ is non-periodic, then $\nu(\theta)$ is a sequence over $\{\0,\1\}$. If $\theta$ is periodic of exact period $n$, then $\nu(\theta)$ also has period $n$ so that the first $n-1$ entries are in $\{\0,\1\}$ and the $n$-th entry is $\*$; such $\nu(\theta)$ are called \emph{$\*$-periodic}. The partition is such that every $\nu(\theta)$ starts with $\1$ (unless $\theta=0$).

Kneading sequences were originally introduced for real quadratic polynomials $x\mapsto x^2+c$ with $c\in\R$. The critical point $0$ separates the dynamical interval into two parts, say $\0$ and $\1$ so that the critical value is in part $\1$. Then the kneading sequence is the sequence of labels $\0$ and $\1$ that the critical orbit visits. The relation to our definition above is as follows: suppose the dynamic ray at angle $R(\theta)$ lands at the critical value. Then the two rays at angles $\theta/2$ and $(\theta+1)/2$ land at the critical point and separate the complex plane into two parts, and the kneading sequence of the critical value along the dynamical interval clearly equals the kneading sequence of the angle $\theta$ as defined above.

\section{Internal Addresses of the Mandelbrot Set}
\label{Sub:IntAddr}

\begin{figure}[htbp]
\setlength{\unitlength}{0.9mm}
\begin{picture}(80,80)
\includegraphics[width=72mm]{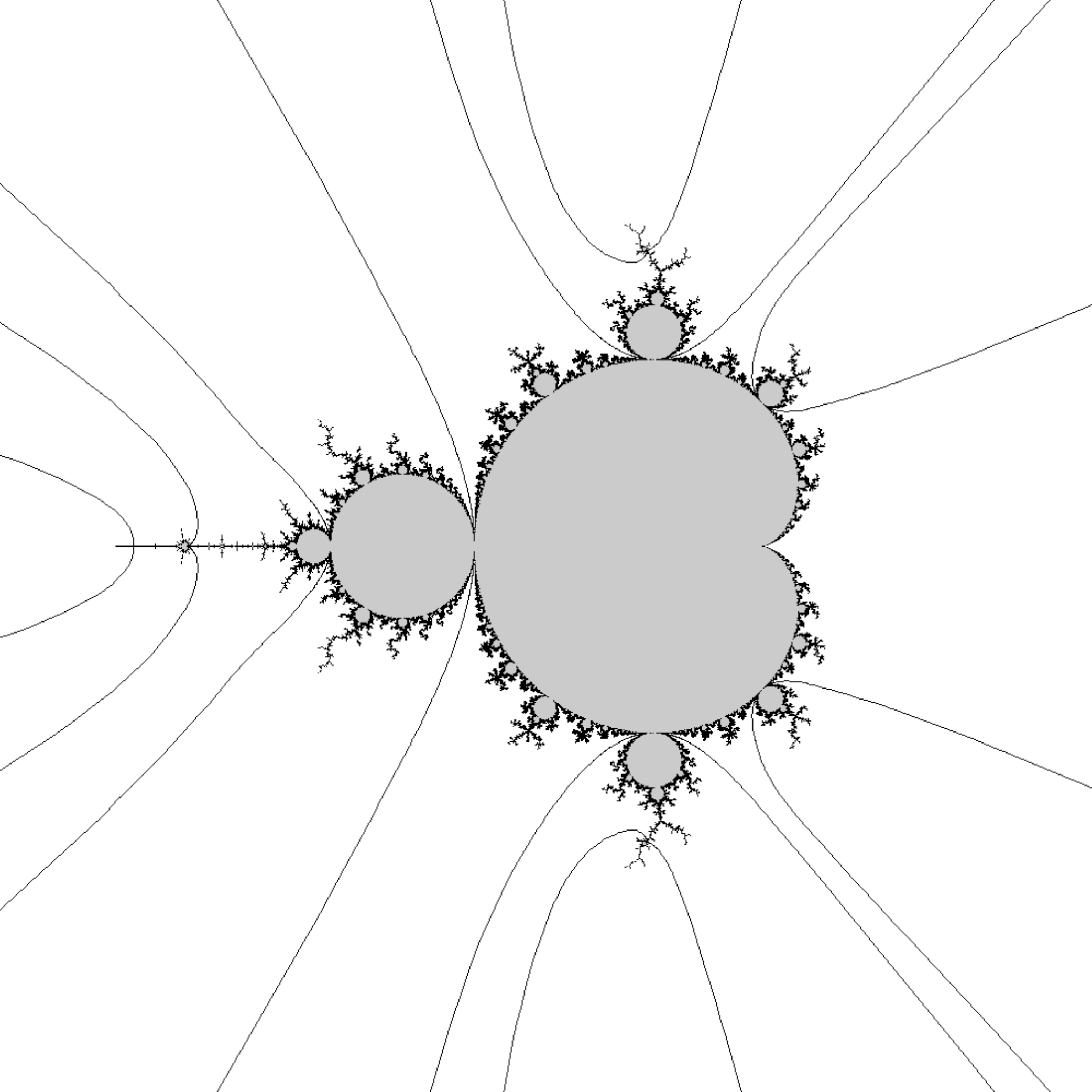} 
\put(-80,42){$7/15$}
\put(-80,50){$3/7$}
\put(-80,60){$6/15$}
\put(-70,77){$1/3$}
\put(-55,77){$2/7$}
\put(-42,77){$4/15$}
\put(-31,77){$3/15$}
\put(-15,77){$1/7$}
\put(-10,72){$2/15$}
\put(-10,53){$1/15$}

\put(-80,36){$8/15$}
\put(-80,28){$4/7$}
\put(-80,18){$9/15$}
\put(-70,0){$2/3$}
\put(-55,0){$5/7$}
\put(-43,0){$11/15$}
\put(-33,5){$12/15$}
\put(-15,0){$6/7$}
\put(-10,05){$13/15$}
\put(-10,24){$14/15$}

\end{picture}
\includegraphics[width=60mm,trim=0 -42 380 580,clip]{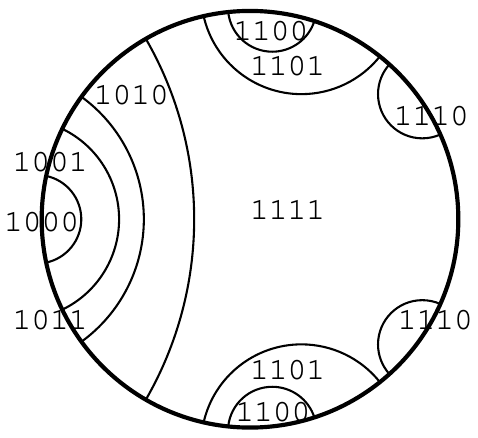}
\caption{Left: the $10$ parameter rays of period up to $4$ in the Mandelbrot set partition $\C$  into $11$ components. In each component, the first $4$ entries of the kneading sequences $\nu(\theta)$ for arbitrary parameter rays $R(\theta)$ are constant. 
Right: the same parameter ray pairs sketched symbolically, and the first $4$ entries in the kneading sequences drawn in.
\hide{\reminder{To be illustrated:} the parameter rays defining the internal address $1\IntAddr 3\IntAddr 4\IntAddr 5\IntAddr 7$ of a hyperbolic component of period $7$.}
}
\label{Fig:PartitionParaRays}
\end{figure}

The parameter rays at periodic angles of periods up to $n$ partition $\C$ into finitely many components.
This partition has interesting symbolic dynamic properties, compare Figure~\ref{Fig:PartitionParaRays} and Lemma~\ref{Lem:ParaRayKneading}: if two parameter rays $R(\theta_1)$ and $R(\theta_2)$ are in the same component, then the kneading sequences $\nu(\theta_1)$ and $\nu(\theta_2)$ associated to $\theta_1$ and $\theta_2$ coincide at least for $n$ entries. In particular, when these two parameter rays land together, then $\theta_1$ and $\theta_2$ must have the same kneading sequence.

The combinatorics of these partitions, and thus the combinatorial structure of the Mandelbrot set, can conveniently be described in terms of \emph{internal addresses}, which are ``human-readable'' recodings of kneading sequences.
In this section we define internal addresses and give their fundamental properties.

\begin{lemma}{Parameter Ray Pairs and Kneading Sequences}
\label{Lem:ParaRayKneading} \lineclear
If two parameter rays $R(\theta)$ and $R(\theta')$ are not
separated by a ray pair of period at most $n$, then $\theta$ and
$\theta'$ have kneading sequences which coincide for at least $n$
entries (provided neither $\theta$ nor $\theta'$ are periodic of
period $n$ or less).

If $R(\theta)$ and $R(\theta')$ with $\theta<\theta'$ form a ray
pair, then $\nu(\theta)=\nu(\theta')=:\nu_\*$. If in addition
both angles are periodic with exact period $n$, then
$\nu_\*$ is $\*$-periodic of period $n$, we have
\[
\lim_{\phi\nearrow\theta}\nu(\phi)=
\lim_{\phi\searrow\theta'}\nu(\phi) \qquad\mbox{ and } \qquad
\lim_{\phi\searrow\theta}\nu(\phi)=
\lim_{\phi\nearrow\theta'}\nu(\phi) \,\,,
\]
and both limits are periodic with period $n$ or dividing $n$ so that their
first difference is exactly at the $n$-th position.
\end{lemma}

\proof
As $\theta$ varies in $\Circle$, the
$n$-th entry in its kneading sequence changes exactly at angles
which are periodic of period dividing $n$.

Consider two external angles $\theta<\theta'$ which are not
separated by any ray pair of period at most $n$, and which are not periodic
of period up to $n$ (it is allowed that $(\theta,\theta')$ contains  periodic angles of period up to $n$, as long as the other angle of the same ray pair is also contained
in $(\theta,\theta')$). Then for every $k\leq n$, the parameter rays of period $k$ with angles
$\phi\in(\theta,\theta')$ must land in pairs 
(Theorem~\ref{Thm:PropertiesParaRayPairs}), so the number of
such rays is even. Therefore, as the angle varies from $\theta$
to $\theta'$, the $k$-th entry in the kneading sequence of
$\theta$ changes an even number of times, and the kneading
sequences of $\theta$ and $\theta'$ have identical $k$-th
entries. This settles the first claim, and it shows that
$\nu(\theta)=\nu(\theta')$ if $R(\theta)$ and $R(\theta')$ land
together and neither $\theta$ nor $\theta'$ are periodic.

However, if $R(\theta)$ and $R(\theta')$ land together and one of the two angles is periodic, 
then both are periodic with the same exact period, say $n$
(Theorem~\ref{Thm:PropertiesParaRayPairs}). 
In this case, the kneading sequences
$\nu(\theta)$ and $\nu(\theta')$ are $\*$-periodic with exact
period $n$, so they coincide as soon as their first $n-1$ entries
coincide; this case is covered by the first claim.

It is quite easy to see that  limits such as $\lim_{\phi\searrow\theta}\nu(\phi)$
exist; moreover, if $\theta$ is non-periodic, then the limit equals $\nu(\theta)$. If
$\theta$ is periodic of exact period $n$, then $\lim_{\phi\searrow\theta}\nu(\phi)$ 
and $\lim_{\phi\searrow\theta}\nu(\phi)$ are periodic of
period $n$ or dividing $n$, they contain no $\*$, and their
first $n-1$ entries coincide with those of $\nu(\theta)$. 
Clearly, $\lim_{\phi\searrow\theta}\nu(\phi)$ and
$\lim_{\phi\nearrow\theta}\nu(\phi)$ differ exactly at positions
which are multiplies of $n$.
Finally, $\lim_{\phi\searrow\theta}\nu(\phi)=
\lim_{\phi\nearrow\theta'}\nu(\phi)$ because for sufficiently
small $\eps$, the parameter rays $R(\theta+\eps)$ and
$R(\theta'-\eps)$ are not separated by a ray pair of period at
most $n$, so the limits must be equal. 
\qed

The following should be taken as an algorithmic definition of
internal addresses in the Mandelbrot set. 

\begin{algo}{Internal Address in Parameter Space}
\label{Alg:IntAddrMandel} \lineclear
Given a parameter $c\in\C$, the internal address
of $c$ in the Mandelbrot set is a strictly increasing finite or infinite sequence
of integers. It is defined as follows:
\begin{description}
\item[seed]
the internal address starts with $S_0=1$ and the ray pair
$\RP(0,1)$;
\item[inductive step]
if $S_0\IntAdr\dots\IntAdr S_k$ is an initial segment of the
internal address of $c$, where $S_k$ corresponds to a ray
pair $\RP(\theta_k,\theta'_k)$ of period $S_k$, then let
$\RP(\theta_{k+1},\theta'_{k+1})$ be the ray pair of least
period which separates $\RP(\theta_k,\theta'_k)$ from $c$
or for which $c\in\RP(\theta_{k+1},\theta'_{k+1})$; let $S_{k+1}$ be
the common period of $\theta_{k+1}$ and $\theta'_{k+1}$.  The internal
address of $c$ then continues as $S_0\IntAdr\dots\IntAdr S_k \IntAdr S_{k+1}$
(with the ray pair $\RP(\theta_{k+1},\theta'_{k+1})$).
\end{description}
This continues for every $k\ge 1$ unless there is a finite $k$
so that $\RP(\theta_k,\theta'_k)$ is not separated from
$c$ by any periodic ray pair (in particular if $c\in\RP(\theta_k,\theta'_k)$).
\end{algo}

\remark
The internal address is only the sequence $S_0\IntAdr\dots\IntAdr S_k \IntAdr \dots$ of periods;
it does not contain the ray pairs used in the construction. The ray pair $\RP(\theta_{k+1},\theta'_{k+1}$) of lowest period is always unique by Lavaurs' Lemma 
\version{\ref{Lem:Lavaurs}}{(see Theorem~\ref{Thm:PropertiesParaRayPairs})}.

The internal address of $c\in\M$ can be viewed as a road map description for the way
from the origin to $c$ in the Mandelbrot set: the way begins at the hyperbolic component of period $1$, and  
at any intermediate place, the internal address
describes the most important landmark on the remaining way to $c$; see Figure~\ref{Fig:PartitionParaRays}. Landmarks are
hyperbolic components (or equivalently the periodic parameter rays landing at their roots,
see \version{Theorem~\ref{Thm:HypComps}),}{Theorem~\ref{Thm:PropertiesParaRayPairs}),} 
and hyperbolic components are the more important the lower their periods are.
The road description starts with the most important landmark: the component of
period $1$, and inductively continues with the period of the component of lowest period
on the remaining way. 

Different hyperbolic components (or combinatorial classes) are not distinguished completely by their internal addresses; the remaining ambiguity has a combinatorial interpretation and will be removed
by \emph{angled internal addresses}: see Theorem~\ref{Thm:CompletenessAngledIntAddr}.

We can give an analogous definition of internal addresses in dynamic planes.
\begin{defi}{Internal Address in Dynamical Planes}
\label{Def:IntAddrJulia} \lineclear
Consider a polynomial $p_c$ for which all periodic dynamic rays land and let $K_c$ be the filled-in Julia set. For a point $z\in K_c$, the internal address of $z$ is defined as follows, in analogy to Algorithm~\ref{Alg:IntAddrMandel}:
\begin{description}
\item[seed]
the internal address starts with $S_0=1$ and the ray pair
$\RP_c(0,1)$;
\item[inductive step]
if $S_0\IntAdr\dots\IntAdr S_k$ is an initial segment of the
internal address of $z$, where $S_k$ corresponds to a dynamic ray
pair $\RP_c(\theta_k,\theta'_k)$ of period $S_k$, then let
$\RP_c(\theta_{k+1},\theta'_{k+1})$ be the dynamic ray pair of least
period which separates $\RP_c(\theta_k,\theta'_k)$ from $z$
or for which $z\in\RP_c(\theta_{k+1},\theta'_{k+1})$; let $S_{k+1}$ be the
common period of $\theta_{k+1}$ and $\theta'_{k+1}$. The internal
address of $z$ then continues as $S_0\IntAdr\dots\IntAdr S_k \IntAdr S_{k+1}$ (with the ray pair $\RP_c(\theta_{k+1},\theta'_{k+1})$).%
\footnote{If several dynamic ray pairs of equal period $S_{k+1}$ separate $P_c(\theta_k,\theta'_k)$ from $z$, take the one which minimizes $|\theta'_{k+1}-\theta_{k+1}|$. (There is an analog to Lavaurs' Lemma~\version{\ref{Lem:Lavaurs}}{} in dynamical planes which says that all candidate ray pairs have to land at the same periodic point; this is not hard to prove, but we do not need it here.)}
\end{description}
This continues for every $k\ge 1$ unless there is a finite $k$
so that $\RP_c(\theta_k,\theta'_k)$ is not separated from
$z$ by any periodic ray pair.
\end{defi}

Every kneading sequence has an associated internal address as follows%
\version{ (compare Definition~\ref{DefKneading})}{}:

\begin{defi}{Internal Address and Kneading Sequence}
\label{Def:IntAddrKneading} \lineclear
Given a kneading sequence $\nu$ with initial entry $\1$, it has the following associated internal address $S_0\IntAddr S_1\IntAddr \dots\IntAddr S_k\IntAddr\dots$: we start with $S_0=1$ and $\nu_0=\ovl{\1}$. Then define recursively $S_{k+1}$ as the position of the first difference between $\nu$ and $\nu_k$, and let $\nu_{S_{k+1}}$ be the periodic continuation of the first $S_{k+1}$ entries in $\nu$ (if $\nu$ is periodic or period $S_k$, then the internal address is finite and stops with entry $S_k$).
\end{defi}
Observe that this definition is algorithmic. It can be inverted in the obvious way so as to assign to any finite or infinite strictly increasing sequence starting with $1$ (viewed as internal address) a kneading sequence consisting of entries $\0$ and $\1$ and starting with $\1$\version{ (Algorithm~\ref{AlgIntAdrKneading})}{}.

For a $\*$-periodic kneading sequence $\nu$ of exact
period $n$, we \version{defined}{define} $\A(\nu)$ and $\Abar(\nu)$ as the two
sequences in which every $\*$ \version{was}{is} replaced consistently by $\0$ or
consistently by $\1$, chosen so that the internal address of
$\A(\nu)$ contains the entry $n$, while the internal address of
$\Abar(\nu)$ does not\version{ (Definition~\ref{DefUpperLower})}{}. 
It turns out that $\A(\nu)$ has exact period $n$, while the exact period of
$\Abar(\nu)$ equals or divides $n$ 
\version{\reminder{Ref to 7.5 of this paper?}
(Proposition~\ref{PropExactPeriod} and independently Lemma~\ref{LemExactPeriod})}
{(see the proof of Lemma~\ref{Lem:PermutationsSymDyn}, \cite[Lemma~4.3]{TreesExist} or 
\cite[{\PropExactPeriod} and independently \LemExactPeriod]{SymDyn})}.
The sequences $\A(\nu)$ and $\Abar(\nu)$ are called the
{\em upper} and {\em lower} periodic kneading sequences associated to $\nu$.

The point of the various algorithmic definitions of internal addresses is of course the following.

\begin{prop}{Equal Internal Addresses}
\label{Prop:EqualIntAddr} \lineclear
(1)
For every $c\in\C$ so that all periodic dynamic rays of $p_c$ land, the internal address in parameter space equals the
internal address of the critical value of $c$ in the dynamical plane of $p_c$.

(2) 
Let $\RP(\theta_k,\theta'_k)$ be the ray pairs from
Algorithm~\ref{Alg:IntAddrMandel}, let $S_k$ be their periods 
and let $\nu^k$ be the common kneading sequence of\/ $\theta_k$ and $\theta'_k$. 
Then each $\nu^k$ is a $\*$-periodic kneading sequence of period $S_k$ so that 
\[
\lim_{\phi\searrow\theta_k}\nu(\phi)=
\lim_{\phi\nearrow\theta'_k}\nu(\phi)=\A(\nu^k)
\qquad \mbox{ and } \qquad
\lim_{\phi\nearrow\theta_k}\nu(\phi)=
\lim_{\phi\searrow\theta'_k}\nu(\phi)=\Abar(\nu^k) \,\,.
\]

(3) The first $S_k$ entries in $\A(\nu^k)$ coincide with those of $\nu^{k+1}$ and, 
if $c\in R(\theta)$, with the first $S_k$ entries of $\nu(\theta)$.

(4) 
For every $\theta\in\Circle$,
the internal address of any parameter $c\in R(\theta)$ in
Algorithm~\ref{Alg:IntAddrMandel} is the same as the internal
address associated to the kneading sequence of\/ $\theta$
from Definition~\ref{Def:IntAddrKneading}.
\end{prop}

\proof
The first statement is simply a reformulation of Corollary~\ref{Cor:RayCorrespondence}
about the correspondence of parameter ray pairs and characteristic dynamic ray pairs:
all we need to observe is that the internal address of the critical value in the Julia set 
uses only characteristic
dynamic ray pairs; this follows directly from the definition of ``characteristic''. 

We prove the second statement by induction, starting with the ray pair $\RP(\theta_0,\theta'_0)$ 
with $\theta_0=0$, $\theta'_0=1$ and $S_0=1$; both angles $\theta_0$ and $\theta'_0$ 
have kneading sequence $\nu^0=\ovl\*$ and $\A(\nu^0)=\ovl\1$.

For the inductive step, suppose that $\RP(\theta_k,\theta'_k)$ is a ray pair of period $S_k$
with $\theta_k<\theta'_k$, $\nu(\theta_k)=\nu(\theta'_k)=\nu^k$ and
$\lim_{\phi\searrow\theta_k}\nu(\phi)=\lim_{\phi\nearrow\theta'_k}\nu(\phi)=\A(\nu^k)$, and
$c$ is not separated from $\RP(\theta_k,\theta'_k)$ by a ray pair of period $S_k$ or less.

As in Algorithm~\ref{Alg:IntAddrMandel}, let $\RP(\theta_{k+1},
\theta'_{k+1})$ be a ray pair of lowest period, say $S_{k+1}$,
which separates $\RP(\theta_k,\theta'_k)$ from
$c$ (or which contains $c$); then $\theta_k<\theta_{k+1}<\theta'_{k+1}<\theta'_k$.
We have $S_{k+1}> S_k$ by construction, 
and the ray pair $\RP(\theta_{k+1},\theta'_{k+1})$ is unique by Lavaurs' Lemma\version{~\ref{Lem:Lavaurs}}{}. 
Since $R(\theta_k)$ and $R(\theta_{k+1})$ are not separated by a ray pair of period
$S_{k+1}$ or less, it follows that the first $S_{k+1}$ entries in
$\lim_{\phi\searrow\theta_k}\nu(\phi)=\A(\nu^k)$ and in
$\lim_{\phi\nearrow \theta_{k+1}}\nu(\phi)$ are equal (the same
holds for $\lim_{\phi\nearrow \theta'_k}\nu(\phi)=\A(\nu^k)$ and
$\lim_{\phi\searrow \theta'_{k+1}}\nu(\phi)$). 

Now we show that
$\lim_{\phi\nearrow\theta_{k+1}}\nu(\phi)
=\lim_{\phi\searrow\theta'_{k+1}}\nu(\phi)=\Abar(\nu^{k+1})$; the
first equality is Lemma~\ref{Lem:ParaRayKneading}.
The internal address associated to $\nu(\theta_k+\eps)$ has no entries in $\{S_k+1,\dots,S_{k+1}\}$ provided $\eps>0$ is sufficiently small (again Lemma~\ref{Lem:ParaRayKneading}). The internal address associated to $\nu(\theta_{k+1}-\eps)$ can then have no entry in $\{S_k+1,\dots,S_{k+1}\}$  for small $\eps$ either (or $R(\theta_k)$ and $R(\theta_{k+1})$ would have to be separated by
a ray pair of period at most $S_{k+1}$). 
Thus $\lim_{\phi\nearrow\theta_{k+1}}\nu(\phi)=\Abar(\nu^{k+1})$.
The other limiting kneading sequences $\lim_{\phi\searrow\theta_{k+1}}\nu(\phi)$ and
$\lim_{\phi\nearrow\theta'_{k+1}}\nu(\phi)$ must then be equal to $\A(\nu^{k+1})$.

Claim (3) follows because neither $\RP(\theta_{k+1},\theta'_{k+1})$ nor $c$ are separated from
$\RP(\theta_k,\theta'_k)$ by a ray pair of period $S_k$ or less.

We prove Claim (4) again by induction. If $c\in R(\theta)$, assume by induction that the internal address 
of $\nu(\theta)$ starts with $1\IntAddr\dots\IntAddr S_k$ and $\theta_k<\theta<\theta'_k$; then $\nu(\theta)$ coincides with $\A(\nu^k)$ for at least $S_k$ entries.
The ray pair of least period separating $R(\theta_k)$ and $R(\theta)$ is $\RP(\theta_{k+1},\theta'_{k+1})$, so the first difference between $\A(\nu^k)$ and $\nu(\theta)$ occurs at position $S_{k+1}$. Hence the internal address of $\nu(\theta)$ as in Definition~\ref{Def:IntAddrKneading} continues as $1\IntAddr\dots\IntAddr S_k\IntAddr S_{k+1}$ and we have  $\theta_{k+1}<\theta<\theta'_{k+1}$ as required to maintain the induction.
\qed

The previous proof also shows the following useful observation.
\begin{coro}{Intermediate Ray Pair of Lowest Period}
\label{Cor:IntermediateRayPair} \lineclear
Let $\RP(\theta_1,\theta'_1)$ and $\RP(\theta_2,\theta'_2)$ be two parameter ray pairs (not necessarily at periodic angles) and suppose that $\RP(\theta_1,\theta'_1)$ separates $\RP(\theta_2,\theta'_2)$ from the origin. If the limits $\lim_{\phi\nearrow\theta_2}\nu(\phi)$ and $\lim_{\phi\searrow\theta_1}\nu(\phi)$ do not differ, then the two ray pairs $\RP(\theta_1,\theta'_1)$ and $\RP(\theta_2,\theta'_2)$ are not separated by any periodic ray pair. If the limits do differ, say at position $n$ for the first time, then both ray pairs are separated by a unique periodic ray pair $\RP(\theta,\theta')$ of period $n$ but not by ray pairs of lower period; the first $n-1$ entries in $\nu(\theta)=\nu(\theta')$ coincide with those of $\lim_{\phi\nearrow\theta_2}\nu(\phi)$ and of $\lim_{\phi\searrow\theta_1}\nu(\phi)$ (while the $n$-th entries are of course $\*$).
\end{coro}
\proof
Let $n$ be the position of the first difference in $\lim_{\phi\searrow\theta_1}\nu(\phi)$ and $\lim_{\phi\nearrow\theta_2}\nu(\phi)$. Then the number of periodic angles of period dividing $n$ in $(\theta_1,\theta_2)$ is odd (because only at these angles, the $n$-th entry in the kneading sequence changes). Since parameter rays at period dividing $n$ land in pairs, there is at least one parameter ray pair at period dividing $n$ which is part of a parameter ray pair separating $\RP(\theta_1,\theta'_1)$ and $\RP(\theta_2,\theta'_2)$. If $\RP(\theta_1,\theta'_1)$ and $\RP(\theta_2,\theta'_2)$ are separated by a parameter ray pair of period less than $n$, let $n'$ be the least such period. By Lavaurs' Lemma\version{~\ref{Lem:Lavaurs}}{}, there would be a single such parameter ray pair of period $n'$, and it follows $n'=n$. The remaining claims follow in a similar way.
\qed

Most internal addresses are infinite; exceptions are related to hyperbolic components as follows.
\begin{coro}{Finite Internal Address}
\label{Cor:FiniteIntAddr} \lineclear
The internal address of $c\in\M$ is finite if and only if there is a hyperbolic component $W$ with $c\in\ovl W$. More precisely, if $c\in\ovl W$ but $c$ is not the root of a hyperbolic component other than $W$, then the internal address of $c$ terminates with the period of $W$; if $c$ is the root of a hyperbolic component $W'\neq W$, then the internal address of $c$ terminates with the period of $W'$.

For a non-periodic external angle $\theta\in\Circle$, the internal address of $R(\theta)$ is finite if and only if $R(\theta)$ lands on the boundary of a hyperbolic component $W$.
\end{coro}
\proof
Consider a parameter ray $R(\phi)$ with finite internal address so that $\phi$ is non-periodic and let $n$ be the last entry in this internal address. Then there is a parameter ray pair $\RP(\theta,\theta')$ of period $n$ which separates $R(\phi)$ from the origin, and no periodic parameter ray pair separates $R(\phi)$ from $\RP(\theta,\theta')$. The ray pair $\RP(\theta,\theta')$ bounds the wake of a hyperbolic component $W$ and lands at the root of $W$, so $R(\phi)$ is in the wake of $W$ but not in one of its subwakes. Every such parameter ray $R(\phi)$ lands on the boundary of $W$; see
\version{Corollary~\ref{Cor:IrrationalIntAngles}}{\cite{HY}\hide{or \cite[\CorIrrationalIntAngles]{SymDyn}}}. The converse is clear.

The statements about $c\in\M$ follow in a similar way, using the fact that the limb of $W$ is the union of $\ovl W$ and its sublimbs at rational internal angles; see \version{Theorem~\ref{Thm:NoGhostLimbs}}{\cite{HY} or \cite[\ThmNoGhostLimbs]{SymDyn}}.
\qed

\remark
The previous result illustrates (through the topology of the Mandelbrot set) the fact that non-periodic external angles can generate periodic kneading sequences which have finite internal addresses. There are periodic kneading sequences which have infinite internal addresses (for instance, $\ovl{\1\0\1}$); but these are never generated by external angles; this essentially follows from the previous result (see also \version{(Lemma~\ref{LemPeriodicItineraryUpper})}{\cite[\LemPeriodicItineraryUpper]{SymDyn}}).

For a hyperbolic component $W$, we have many different ways of associating an internal address to it: \begin{enumerate}
\item
using Algorithm~\ref{Alg:IntAddrMandel} within the complex plane of the Mandelbrot set; 
\item
taking a parameter ray $R(\theta)$ at a periodic angle that lands at the root of $W$ and then Algorithm~\ref{Alg:IntAddrMandel} for this parameter ray;
\item
taking a parameter ray $R(\theta')$ at an irrational angle that lands at a  point in $\partial W$ (as in Corollary~\ref{Cor:FiniteIntAddr}) and then again Algorithm~\ref{Alg:IntAddrMandel} for this parameter ray;
\item
taking a parameter ray $R(\theta)$ or $R(\theta')$ as in (2) or (3) and then using the internal address associated to the kneading sequence of $\theta$ or $\theta'$;  
\version{see Definitions~\ref{DefKneading} and \ref{Def:IntAddrKneading}}
{see Definition~\ref{Def:IntAddrKneading}}. 
The periodic angle $\theta$ generates a $\*$-periodic kneading sequence, while the non-periodic angle $\theta'$ generates a periodic kneading sequence without $\*$;
\item
for every $c\in W$, we have the internal address of the critical value in the Julia set from Definition~\ref{Def:IntAddrJulia};
\item
specifically if $c$ is the center of $W$, then the critical orbit is periodic and $p_c$ has a Hubbard tree in the original sense of Douady and Hubbard 
\version{(Definition~\ref{DefDouadyHubbardTree})}{\cite{Orsay}}, so we can use all definitions from 
\version{Proposition~\ref{PropIntAddr}}{\cite[\PropIntAddr]{SymDyn}}.
\end{enumerate}
All these internal addresses are of course the same: this is obvious for the first three; the next two definitions agree with the first ones by Proposition~\ref{Prop:EqualIntAddr}~(4) and (1), and the last two agree by \version{Proposition~\ref{PropIntAddr}}{\cite[\PropIntAddr]{SymDyn}}.


\begin{defi}{Angled Internal Address}
\label{Def:AngledIntAddr} \lineclear
For a parameter $c\in\C$, the {\em angled internal address} is the sequence
\[
(S_0)_{p_0/q_0}\IntAddr (S_1)_{p_1/q_1}\IntAddr \dots \IntAddr (S_k)_{p_k/q_k}\IntAddr\dots
\]
where $S_0\IntAddr S_1\IntAddr\dots\IntAddr S_k\IntAddr\dots$ is the internal address
of $c$ as in Algorithm~\ref{Alg:IntAddrMandel} and the angles $p_k/q_k$ are defined as follows:
for $k\ge 0$, let $\RP(\theta_k,\theta'_k)$ be the parameter ray pair associated to the
entry $S_k$ in the internal address of $c$; then the landing point of $\RP(\theta_k,\theta_k')$ 
is the root of a hyperbolic
component $W_k$ of period $S_k$. The angle $p_k/q_k$ is defined such that
$c$ is contained in the closure of the $p_k/q_k$-subwake of $W_k$.

If the internal address of $c$ is finite and terminates with an entry $S_k$ (which happens if and only if $c$ is not contained in the closure of any subwake of $W_k$), then
the angled internal address of $c$ is also finite and terminates with the same entry $S_k$
without angle:
$(S_0)_{p_0/q_0}\IntAddr (S_1)_{p_1/q_1}\IntAddr \dots \IntAddr (S_{k-1})_{p_{k-1}/q_{k-1}}\IntAddr (S_k)$.
\end{defi}

\begin{figure}[htbp]
\centerline{\framebox{\includegraphics[scale=1.1]{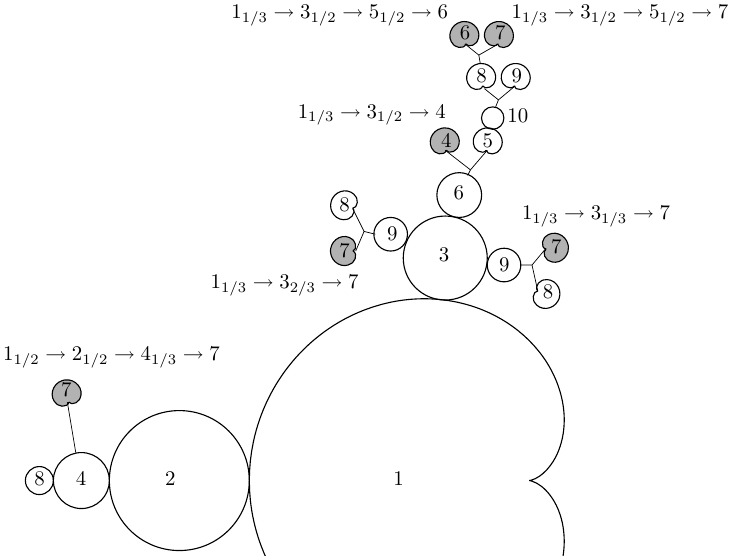}}}
\caption{Angled internal addresses for various hyperbolic components in $\M$.}
\label{Fig:AngledIntAddr}
\end{figure}

This definition is illustrated in Figure~\ref{Fig:AngledIntAddr}.
The main point in this definition is that it distinguishes different points in the Mandelbrot set. A precise statement is given in terms  combinatorial classes:
a \emph{combinatorial class} is the equivalence class of parameters in $\M$ so that two parameters $c_1$ and $c_2$ are equivalent if and only if for both polynomials, the same periodic and preperiodic dynamic land at common points. Equivalently, it is a maximal subset of $\M\sm\{\mbox{parabolic parameters}\}$ so that no two of its points can be separated by a parameter ray pair at periodic angles; a parabolic parameter $c_0$ belongs to the combinatorial class of the hyperbolic component of which $c_0$ is the root; see for instance 
 \version{(see Section~\ref{Sub:CombClass} or \cite[Section~8]{Fibers2})}{\cite[Section~8]{Fibers2}}). Local connectivity of $\M$ is equivalent to the conjecture that all non-hyperbolic combinatorial classes are points.

We can now describe precisely which points in $\M$ angled internal addresses distinguish.

\begin{theo}{Completeness of Angled Internal Addresses}
\label{Thm:CompletenessAngledIntAddr} \lineclear
Two parameters in $\M$ have the same angled internal address if an only if they belong to the same combinatorial class.
In particular, two hyperbolic parameters in $\M$ have the same angled internal address if and only if they belong to the same hyperbolic component.
\end{theo}

We postpone the proof of this theorem and of subsequent results%
\version{}{ to Section~\ref{Sub:ProofsIntAddr}} and first discuss some interesting consequences.

\version
{\heading{The Geometry of Internal Addresses}}
{\subsection{The Geometry of Internal Addresses}}
\label{Sub:GeometryIntAddr}

Internal addresses give an efficient and convenient way to locate and describe parameters in $\M$. For instance, they make it possible to give an answer to the folklore question how to tell where in the Mandelbrot set a given parameter ray lands.

It turns out that the internal address without angles completely determines the denominators of any associated angled internal address.

\begin{lemma}{Estimates on Denominators in Internal Address}
\label{Lem:EstimatesDenominators} \lineclear
In any internal address $(S_0)_{p_0/q_0}\to \dots\to (S_k)_{p_k/q_k}\to (S_{k+1})_{p_{k+1}/q_{q+1}}\to\dots$, 
the denominators $q_k$ satisfy $S_{k+1}/S_k\le q_k<S_{k+1}/S_k+2$; moreover, $q_k=S_{k+1}/S_k$ if and only if the latter is an integer. 

If $S_{k+1}$ is a multiple of $S_k$, then the component of period $S_{k+1}$ is a bifurcation from that of period $S_k$.
\end{lemma}

This lemma provides a simple estimate on the denominators that allows only two possibilities for all $q_k$. 

Here we give a precise combinatorial formula for $q_k$.
\version{ It is analogous to Formula (\ref{numberarms}) in Proposition~\ref{PropArmsOri}.}{}%
Recall that every internal address has an associated kneading sequence $\nu=\nu_1\nu_2\nu_3\dots$ (see Definition~\ref{Def:IntAddrKneading}\version{ or Algorithm~\ref{AlgIntAdrKneading}}{}; this does not depend on the angles). 
To this kneading sequence $\nu$ we define an associated function $\rho_\nu$ as follows: for $r\ge 1$, let 
\[
\rho_\nu(r):=\min\{k>r\colon  \nu_{k}\neq\nu_{k-r}\} \,\,.
\]
The $\rho_\nu$-orbit of $r$ is denoted $\orb_{\rho_\nu}(r)=\{r,\rho_\nu(r),\rho_\nu(\rho_\nu(r)),\dots\}$%
\version{ (compare Definition~\ref{DefRho})}{}. We often write $\rho$ for $\rho_\nu$.

\begin{lemma}{Denominators in Angled Internal Address}
\label{Lem:Denominators} \lineclear
In an angled internal address $(S_0)_{p_0/q_0}\IntAddr \dots \IntAddr (S_k)_{p_k/q_k}\IntAddr (S_{k+1})_{p_{k+1}/q_{k+1}}\dots$, the denominator $q_k$ in the bifurcation angle is uniquely determined by the internal address $S_0\IntAddr\dots\IntAddr S_k\IntAddr S_{k+1}\dots$ as follows: let $\nu$ be the kneading sequence associated to the internal address  and let $\rho$ be the associated function as just described. Let $r\in\{1,2,\ldots,S_k\}$ be congruent to $S_{k+1}$ modulo $S_k$. Then
\[
q_k := \left\{ \begin{array}{ll}
\frac{S_{k+1}-r}{S_k} + 1 & \mbox{ if } S_k \in \orb_{\rho}(r)\,\,,\\
\frac{S_{k+1}-r}{S_k} + 2 & \mbox{ if } S_k \notin \orb_{\rho}(r)\,\,.
\end{array} \right.
\]
\end{lemma}

While the internal address completely specifies the denominators in any corresponding angled internal address, it says says nothing about the numerators: of course, not all internal addresses are realized in the Mandelbrot set (not all are (complex) \emph{admissible}; see
\version{Sections~\ref{SecAdmissCondition} and \ref{SecExtAnglesAdmiss}}{\cite{AdmissKneading} or \cite[\SecAdmissTwice]{SymDyn}}), but this is independent of the numerators. 

\begin{theo}{Independence of Numerators in Angled Internal Address}
\label{Thm:NumeratorsArbitrary} \lineclear
If an angled internal address describes a point in the Mandelbrot set, then the numerators $p_k$ can be changed arbitrarily (coprime to $q_k$) and the modified angled internal address still describes a point in the Mandelbrot set.
\end{theo}

In other words, for every hyperbolic component there is a natural bijection between combinatorial classes of the $p/q$-sublimb and $p'/q$-sublimb, for every $q\ge 2$ and all $p,p'$ coprime to $q$. We thus conjectured in 1994 that these two limbs were homeomorphic by a homeomorphism preserving periods of hyperbolic components (but not the embedding into the plane). This conjecture was established recently \cite{DimaHomeos}, based on work of Dudko \cite{DimaDecoThm}. These limbs have been known to be homeomorphic by work of Branner and Fagella~\cite{BF}, but their homeomorphisms (like all constructed by quasiconformal surgery, for instance as in \cite{Riedl}) preserve the embedding into the complex plane and not periods of hyperbolic components.

Here is a way to find the numerators $p_k$ of the internal address of an external angle.
\begin{lemma}{Numerators in Angled Internal Address}
\label{Lem:Numerators} \lineclear
Suppose the external angle $\theta$ has angled internal address $(S_0)_{p_0/q_0}\IntAddr \dots \IntAddr (S_k)_{p_k/q_k}\IntAddr (S_{k+1})_{p_{k+1}/q_{k+1}}\dots$. In order to find the numerator $p_k$, consider the $q_k-1$ angles $\theta$, $2^{S_k}\theta$, $2^{2S_k}\theta$, \dots, $2^{(q_k-2)S_k}\theta$. Then $p_k$ is the number of these angles that do not exceed $\theta$.
\end{lemma}

If $\theta$ is periodic, then it is the angle of one of the two parameter rays landing at the root of a hyperbolic component. Here is a way to tell which of the two rays it is.

\begin{lemma}{Left or Right Ray}
\label{Lem:LeftRightRay} \lineclear
Let $R(\theta)$ and $R(\theta')$ be the two parameter rays landing at the root of a hyperbolic component $W$ of period $n\ge 2$, and suppose that $\theta<\theta'$. Let $b$ and $b'$ be the $n$-th entries in the binary expansions of $\theta$ and $\theta'$. Then 
\begin{itemize}
\item
if the kneading sequence of $W$ has $n$-th entry $\0$, then $b=1$ and $b'=0$;
\item
if the kneading sequence of $W$ has $n$-th entry $\1$, then $b=0$ and $b'=1$.
\end{itemize}
\end{lemma}

For example, the hyperbolic component with internal address $1_{1/3}\IntAddr 3_{1/2}\IntAddr 4$ has kneading sequence $\ovl{\1\1\0\0}$, and the parameter ray $R(4/15)$ lands at its root. The binary expansion of $4/15$ is $0.\ovl{0100}$. The $4$-th entries in the kneading sequence and in the binary expansion of $4/15$ are $\0$ and $0$, so the second ray landing together with $R(4/15)$ has angle $\theta<4/15$ (indeed, the second ray is $R(3/15)$).

\version{}{
We define the \emph{width of the wake} of a hyperbolic component $W$ as $|W|:=|\theta'-\theta|$, where $\RP(\theta,\theta')$ is the parameter ray pair landing at the root of $W$. If $W$ has period $n$, then the width of the $p/q$-subwake of $W$ equals 
\begin{equation}
|W_{p/q}| = |W|(2^n-1)^2/(2^{qn}-1)\,\,;
\label{Eq:WidthWake}
\end{equation}
this is a folklore result, related to Douady's Tuning Formula~\cite{DoAngles}; a proof can be found in \cite[\PropWidthWake]{SymDyn}.} 

The following result complements the interpretation of internal addresses as road descriptions by saying that whenever the path from the origin to a parameter $c\in\M$ branches off from the main road, an entry in the internal address is generated: the way to most parameters $c\in\M$ traverses infinitely many hyperbolic components, but most of them are traversed ``straight on'' and left into the $1/2$-limb.

\begin{prop}{Sublimbs Other Than $1/2$ in Internal Address}
\label{Prop:OtherThanOneHalf} \lineclear
If a parameter $c\in\C$ is contained in the subwake of a hyperbolic component $W$ at internal angle other than $1/2$, then $W$ occurs in the internal address of $c$.  More precisely, the period of $W$ occurs in the internal address, and the truncation of the angled internal address of $c$ up to this period describes exactly the component $W$.
\end{prop}
\proof
Let $n$ be the period of $W$. If $W$ does not occur in the internal address of $c$, then the internal address of $c$ must have an entry $n'<n$ corresponding to a hyperbolic component $W'$ in the wake of $W$. Denoting the width of $W'$ by $|W'|$, we have $|W'|\ge 1/(2^{n'}-1)$.
\version{
By Proposition~\ref{Prop:WidthWake}, if $|W|$ denotes the width of a hyperbolic component $W$ of period $n$, the width of its $p/q$-subwake is  $|W|(2^n-1)^2/(2^{qn}-1)$,}
{By (\ref{Eq:WidthWake}), the width of the $p/q$-subwake of $W$ is  $|W|(2^n-1)^2/(2^{qn}-1)$,}
so we must have 
\[
1/(2^{n'}-1) \le |W'|\le  |W|(2^n-1)^2/(2^{qn}-1) < (2^n-1)^2/(2^{qn}-1)
\]
 or 
$2^{qn}-1 < (2^{n'}-1)(2^n-1)^2 < (2^n-1)^3 < 2^{3n}-1$, hence $q<3$.
\qed

The internal address of a parameter $c$ also tells whether or not this parameter is renormalizable%
\version{ in the sense of Section~\ref{Sec:Renormalization}}{}:
\begin{prop}{Internal Address and Renormalization}
\label{Prop:IntAddrRenorm} \lineclear
Let $c\in\M$ be a parameter with internal address 
$(S_0)_{p_0/q_0}\IntAddr \dots \IntAddr (S_k)_{p_k/q_k}\dots$. 
\begin{itemize}
\item
$c$ is {\em simply renormalizable} of period $n$ if and only if there is a $k$ with $S_k=n$ and $n|S_{k'}$ for $k'\ge k$;
\item
$c$ is {\em crossed renormalizable} of period $n$ if and only if there is an $m$ strictly dividing $n$ so that $S_k=m$ for an appropriate $k$ and all $S_{k'}$ with $k'>k$ are proper multiples of $n$ (in particular, $n$ does not occur in the internal address). In this case, the crossing point of the little Julia sets has period $m$.
\end{itemize}
\end{prop}

\looseness -1
\remark
We can also describe combinatorially whether any given parameter rays $R(\theta)$ and $R(\theta')$ can land at the same point in $\M$: a necessary condition is that they have the same angled internal address. This condition is also sufficient from a combinatorial point of view: suppose $R(\theta)$ and $R(\theta')$ have the same angled internal address. Then both rays accumulate at the same combinatorial class. If the internal address is finite, then $R(\theta)$ and $R(\theta')$ land on the boundary of the same hyperbolic component. Otherwise, both rays accumulate at the same combinatorial class. As soon as it is known that this combinatorial class consists of a single point (which would imply local connectivity of $\M$ at that point \cite{Fibers2}), then both rays land together.

Here is one example: where in the Mandelbrot set does the parameter ray at angle $\theta=22/(2^7-1)$ land? The angle has period $7$, so it must land at the root of a hyperbolic component of period $7$; but there are $63$ such components; which one is it? And which of the two rays landing at this root is it?

The angle has associated kneading sequence $\nu(\theta)=\ovl{\1\1\0\1\0\1\*}$, also of period $7$ (see Definition~\ref{Def:Kneading}). This binary sequences has an associated internal address (given by the simple algorithm in Definition~\ref{Def:IntAddrKneading}): it equals $1\IntAddr 3\IntAddr 5\IntAddr  7$, and it describes where the ray lands (see Figure~\ref{Fig:AngledIntAddr}): start with the unique hyperbolic component at period $1$; the most important landmark from here to the destination point is the component of period $3$ that bifurcates from the period $1$ component. From here, pass to the unique period $5$ component without passing by a component of period $4$ (or the $4$ would appear in the internal address); and there is a unique such component. The wake of this period $5$ component contains a unique period $6$ component and two period $7$ components: one inside and one outside of the wake of a period $6$ component. Our component is the one outside of the period $6$ wake, and the parameter ray at angle $\theta$ lands at its root. 

There are two components of period $3$ bifurcating from the period $1$ component, and we need the one at internal angle $1/3$: this is resolved by \emph{angled} internal addresses. The denominators are determined by Lemma~\ref{Lem:Denominators}, and the numerators by Lemma~\ref{Lem:Numerators}:  the angled internal address is $1_{1/3}\IntAddr 3_{1/2} \IntAddr 5_{1/2}\IntAddr7$.

Finally, there are two rays landing at the root of the period $7$ component, and our ray is the greater of the two (Lemma~\ref{Lem:LeftRightRay}). Internal addresses also specify which components bifurcate from each other, and which ones are renormalizable (Proposition~\ref{Prop:IntAddrRenorm}). 

\version
{\heading{Proofs about Internal Addresses}}
{\subsection{Proofs about Internal Addresses}}
\label{Sub:ProofsIntAddr}

Now we give the proofs of the results stated so far; most of them go back to \cite{IntAdr}. Many of the proofs are based on the concept of {\em long internal addresses}, which show that even though internal addresses themselves are a compact road description, they encode refinements to arbitrary detail. Parameter ray pairs have a partial order as follows:
$\RP(\theta_1,\theta'_1)<\RP(\theta_2,\theta'_2)$ iff
$\RP(\theta_1,\theta'_1)$ separates $\RP(\theta_2,\theta'_2)$ from the origin $c=0$, 
or equivalently from $c=1/4$, the landing point of the ray pair
$\RP(0,1)$. It will be convenient to say $\RP(0,1)<\RP(\theta,\theta')$ for every ray pair
$\RP(\theta,\theta')\neq \RP(0,1)$. For every parameter $c\in\C$, the set of parameter ray pairs
separating $c$ from the origin is totally ordered (by definition, this set always includes the pair
$\RP(0,1)$ as its minimum). 
A similar partial order can be defined for dynamic ray pairs for every $c\in\M$.

\begin{defi}{Long Internal Address in the Mandelbrot Set}
\label{Def:LongIntAddrMandel} \lineclear
For a parameter ray $c\in\C$, consider the set of periodic parameter ray pairs
which separate $c$ from the origin, totally ordered as described above; the 
\emph{long internal address} of $c$ is the collection of periods of these ray pairs
with the same total order.
\end{defi}
The long internal address always starts with the entry $1$, and it is usually infinite and not well-ordered. For most $c\in\C$, many periods appear several times on the long internal address (which means that several ray pairs of equal period separate $c$ from the origin; only the ray pair of least period is unique).

One useful feature of internal addresses is that they completely encode the associated long internal addresses. The following proposition is in fact algorithmic, as its proof shows.

\begin{prop}{Long Internal Address Encoded in Internal Address}
\label{Prop:LongIntAddrEncoded} \lineclear
Any internal address completely encodes the associated long internal address.
\end{prop}
\proof
The internal address is a strictly increasing (finite or infinite) sequence of integers, each of which comes with an associated $\*$-periodic kneading sequence. If it is the internal address of some $c\in\C$, then each entry in the internal address represents a parameter ray pair with this period. Corollary~\ref{Cor:IntermediateRayPair} describes the least period of a periodic parameter ray pair which separates any two given parameter ray pairs, and it also describes the kneading sequence of the ray pair of least period. This allows us to inductively find the periods of all parameter ray pairs of given maximal periods which separate $c$ from the origin, together with the order and kneading sequences of these ray pairs, and in the limit determines the long internal address.
\qed

\remark
Of course, it makes perfect sense to speak of an {\em angled long internal address}, which is a long internal address in which all entries (except the last, if there is a last entry) are decorated with the bifurcation angles  of the sublimb containing the parameter associated with this address. All entries with angles different from $1/2$ are already contained in the ``short'' internal address by Proposition~\ref{Prop:OtherThanOneHalf}, so the angled internal address completely encodes the angled long internal address.

\proofof{Theorem~\ref{Thm:CompletenessAngledIntAddr} (Completeness of angled internal addresses)}  
If two parameters belong to the same combinatorial class, then by definition they cannot be separated by a parameter ray pair at periodic angles, and then their angled internal address is by definition the same.

For the converse, suppose two parameters $c$ and $c'$ have the same angled internal address. If the two sequences of hyperbolic components in these two angled internal addresses coincide, then it easily follows that $c$ and $c'$ are in the same combinatorial class. It thus suffices to prove the claim for periodic ray pairs and thus for hyperbolic components. If the claim is false, then there is a least period $S_k$ for which there are two different hyperbolic components $W$ and $W'$ of period $S_k$ that share the same angled internal address $(S_1)_{p_1/q_1}\IntAddr \dots \IntAddr (S_{k-1})_{p_{k-1}/q_{k-1}}\IntAddr S_k$. By minimality of $S_k$, the ray pair of period $S_{k-1}$ is the same in both internal addresses.

By the Branch Theorem~\version{\ref{Thm:MandelBranch}}{(see \cite{Orsay} or \cite[Theorem~3.1]{Fibers2})}, there are three possibilities: either (1) $W$ is contained in the wake of $W'$ (or vice versa), or (2) there is a hyperbolic component $W_*$ so that $W$ and $W'$ are in two different of its sublimbs, or (3) there is a Misiurewicz-Thurston point $c_*$ so that $W$ and $W'$ are in two different of its sublimbs. (A \emph{Misiurewicz-Thurston} point is a parameter $c_*$ for which the critical orbit is strictly preperiodic; such a point is the landing point of a finite positive number of parameter rays at preperiodic angles \cite{Orsay,ExtRayMandel}, and a \emph{sublimb} of $c_*$ is a component in the complement of these rays and their landing point, other than the component containing $0$.)

(1) In the first case,  there must by a parameter ray pair of period less than $S_k$ separating $W$ and $W'$ by Lavaurs' Lemma\version{~\ref{Lem:Lavaurs}}{}. This ray pair would have to occur in the internal address of $W'$ (between entries $S_{k-1}$ and $S_k$ or perhaps before $S_{k-1}$) but not of $W$, and this is a contradiction.

(2) The second possibility is handled by angled internal addresses: at least one of $W$ or $W'$ must be contained in a sublimb of $W_*$ other than the $1/2$-sublimb, so by Proposition~\ref{Prop:OtherThanOneHalf}, $W_*$ must occur in the internal address, and the angles at $W_*$ in the angled internal address distinguish $W$ and $W'$.

(3) The case of a Misiurewicz-Thurston $c_*$ point needs more attention. If  $c_*$ has at least two sublimbs, then it is the landing point of $k\ge 3$ parameter rays $R(\theta_1)$ \dots $R(\theta_k)$ at preperiodic angles.

By Proposition~\ref{Prop:LongIntAddrEncoded}, $W$ and $W'$ have the same long internal address, and by minimality of $S_k$ there is no hyperbolic component $W^*$ of period less than $S_k$ in a sublimb of $c^*$ so that $W\subset \wake{W^*}$. 

Let $\RP(\theta,\theta')$ be the parameter ray pair landing at the root of  $W$ and let $c$ be the center of $W$. In the dynamics of $p_c$, there is a characteristic preperiodic point $w$ which is the landing point of the dynamic rays $R_c(\theta_1)$, \dots $R_c(\theta_K)$ 
(the definition of characteristic preperiodic points is in analogy to the definition of characteristic periodic points after Theorem~\ref{Thm:PropertiesParaRayPairs}; the existence of the characteristic preperiodic point $w$ is well known (``The Correspondence Theorem''); see e.g., \version{Theorem~\ref{Thm:Correspondence}}{\cite[Theorem~2.1]{Fibers2} or \cite[\ThmCorrespondencence]{SymDyn}}). We use the Hubbard tree $T_c$ of $p_c$ in the original sense of 
Douady and Hubbard 
\version{(Definition~\ref{DefDouadyHubbardTree})}{\cite{Orsay}}. Let $I\subset T_c$ be the arc connecting $w$ to $c$. 

If the restriction $p_c^{\circ S_k}|_I$ is not injective, then let $n\le S_k$ be maximal so that $p_c^{\circ (n-1)}|_I$ is injective. Then there is a sub-arc $I'\subset I$ starting at $w$ so that $p_c^{\circ n}|_{I'}$ is injective and $p_c^{\circ n}(I')$ connects $p_c^{\circ n}(w)$ with $c$. If $n=S_k$ then $I'=I$ because $c$ is an endpoint of $p_c^{\circ n}(I)$, a contradiction; thus $n<S_k$. Since $w$ is characteristic, this implies that $p_c^{\circ n}(I')\supset I'$, so there is a fixed point $z$ of $p_c^{\circ n}$ on $I'$. If $z$ is not characteristic, then the characteristic point on the orbit of $z$ is between $z$ and $c$. In any case, we have a hyperbolic component $W^*$ of period $n<S_k$ in a sublimb of $c^*$ so that $W\subset \wake{W^*}$ (again by Theorem~\ref{Thm:PropertiesParaRayPairs}), and this is a contradiction.

It follows that the restriction $p_c^{\circ S_k}|_I$ is injective. There is a unique component of $\C\sm \left(R_c(2^{S_k}\theta_1)\cup\dots\cup R_c(2^{S_k}\theta_k)\cup\{p_c^{\circ S_k}(w)\}\right)$ that contains $p_c^{\circ S_k}(I)$: it is the component containing $c$ and the dynamic ray pair $\RP_c(\theta,\theta')$, and thus also the dynamic rays $R_c(\theta_1),\dots,R_c(\theta_k)$, so this component is uniquely specified by the external angles of $c_*$ together with $S_k$. But by injectivity of $p_c^{\circ S_k}|_I$, this also uniquely specifies the component of $\C\sm (R_c(\theta_1)\cup\dots\cup R_c(\theta_k)\cup\{z\})$ that must contain $I$ and hence $c$. In other words, the subwake of $c_*$ containing $W$ (and, by symmetry, $W'$) is uniquely specified.
\qed

\proofof{Lemma~\ref{Lem:EstimatesDenominators} (Estimates on Denominators)}
The angled internal address $(S_0)_{p_0/q_0}\IntAddr (S_1)_{p_1/q_1}\IntAddr \dots \IntAddr (S_k)_{p_k/q_k}\IntAddr (S_{k+1})_{p_{k+1}/q_{k+1}}\dots$, when truncated at periods $S_k$, describes a sequence of hyperbolic components $W_k$ of periods $S_k$. If $W_{k+1}$ is contained in the $p_k/q_k$-sublimb of $W_k$, then we have $S_{k+1}\le q_kS_k$ (otherwise an entry $q_kS_k$ would be generated in the internal address); thus $q_k\ge S_{k+1}/S_k$. The other inequality follows from \version{Proposition~\ref{Prop:WidthWake}:}{(\ref{Eq:WidthWake}):}
the width of the wake of $W_{k+1}$ cannot exceed the width of the $p_k/q_k$-subwake of $W_k$, so $1/(2^{S_{k+1}}-1)\le (2^{S_k}-1)^2/(2^{q_kS_k}-1)$ or 
\[
2^{q_kS_k}-1 \le (2^{S_{k+1}}-1) (2^{S_k}-1)^2 < 2^{S_{k+1}+2S_k}-1 \,\,,
\]
hence $S_{k+1}>(q_k-2)S_k$ or $q_k<S_{k+1}/S_k+2$.

It remains to show that whenever $S_{k+1}/S_k$ is an integer, it equals $q_k$: there are associated parameter ray pairs $\RP(\theta_k,\theta'_k)$ of period $S_k$ and $\RP(\theta_{k+1},\theta'_{k+1})$ of period $S_{k+1}$, and the limiting kneading sequences $\lim_{\phi\searrow\theta_k}\nu(\phi)$ and $\lim_{\phi\nearrow\theta_{k+1}}\nu(\phi)$ are periodic of period $S_k$ and $S_{k+1}$; both can be viewed as being periodic of period $S_{k+1}$. Since they correspond to adjacent entries in the internal address, these ray pairs cannot be separated by a ray pair of period up to $S_{k+1}$, so both limiting kneading sequences are equal. Therefore, Corollary~\ref{Cor:IntermediateRayPair} implies that the two ray pairs are not separated by any periodic parameter ray pair. If $W_{k+1}$ was not a bifurcation from $W_k$, then it would have to be in some subwake of $W_k$ whose boundary would be some ray pair separating $\RP(\theta_k,\theta'_k)$ from $\RP(\theta_{k+1},\theta'_{k+1})$; this is not the case. So let $p_k/q_k$ be the bifurcation angle from $W_k$ to $W_{k+1}$; then the corresponding periods satisfy $S_{k+1}=q_kS_k$ as claimed.
\qed

In the following lemma, we use the function $\rho_\nu$ as defined before Lemma~\ref{Lem:Denominators}. 

\begin{lemma}{Intermediate Ray Pair of Lowest Period}
\label{Lem:IntermediateRayPairSpecialCase} \lineclear
Let $\RP(\theta_k,\theta_k')$ and $\RP(\theta_{k+1},\theta'_{k+1})$ be two periodic parameter ray pairs with periods $S_k<S_{k+1}$ and suppose that $\RP(\theta_k,\theta_k')$ separates $\RP(\theta_{k+1},\theta'_{k+1})$ from the origin. Write $S_{k+1}=aS_k+r$ with $r\in\{1,\dots,S_k\}$. Let $S$ be the least period of a ray pair separating $\RP(\theta_k,\theta_k')$ from $\RP(\theta_{k+1},\theta'_{k+1})$. If $S_{k+1}<S\le(a+1)S_k$, then $S=aS_k+\rho_\nu(r)$, where $\nu$ is any kneading sequence that has the same initial $S_k$ entries as $\nu(\theta_{k+1})=\nu(\theta'_{k+1})$.
\end{lemma}
\proof
Let $B$ and $R$ denote the initial blocks in $\A(\nu(\theta_k))=\A(\nu(\theta'_k))$ consisting of the first $S_k$ or $r$ entries, respectively. Then $\A(\nu(\theta_k))=\ovl B$. Since $S>S_{k+1}$, Corollary~\ref{Cor:IntermediateRayPair} implies $\Abar(\nu(\theta_{k+1}))= \ovl{B^aR}$. Therefore,  $S$ is the position of the first difference between  $\ovl B$ and $\ovl{B^aR}$, and $S-aS_k$ is the position of the first difference between $\ovl B$ and $\ovl{RB^a}$. Since $S\le (a+1)S_k$, this difference occurs among the first $S_k$ symbols, and these are specified by $B$ and $RB$. Since $R$ is the initial segment of $B$ of length $r$, $S-aS_k$ equals $\rho_{\nu}(r)$ for any sequence $\nu$ that starts with $B$.
\qed

\proofof{Lemma~\ref{Lem:Denominators} (Finding denominators)}
The internal address (without angles) uniquely determines the long internal address by Proposition~\ref{Prop:LongIntAddrEncoded}. The entry $S_k$ occurs in the internal address and the subsequent entry in the long internal address is $q_kS_k$, so the denominators are uniquely encoded (and depend only on the $q_kS_k$ initial entries in the kneading sequence). Recall the bound  $S_{k+1}\le q_kS_k < S_{k+1}+2S_k$ from Lemma~\ref{Lem:EstimatesDenominators}.

Write again $S_{k+1}=aS_k+r$ with $r\in\{1,\dots,S_k\}$. If $r=S_k$, then $S_{k+1}$ is divisible by $S_k$ and $q_k=S_{k+1}/S_k$ by Lemma~\ref{Lem:EstimatesDenominators}, and this is what our formula predicts. Otherwise, we have $q_k\in\{a+1,a+2\}$. Below, we will find the lowest period $S'$ between $S_k$ and $S_{k+1}$, then the lowest period $S''$ between $S_k$ and $S'$, and so on
(of course, the ``between'' refers to the order of the associated ray pairs). This procedure must eventually reach the bifurcating period $q_kS_k$. If $q_kS_k=(a+1)S_k$, then eventually one of the periods $S'$, $S''$,\dots must be equal to $(a+1)S_k$. If not, the sequence $S'$, $S''$, \dots skips $(a+1)S_k$, and then necessarily $q_kS_k=(a+2)S_k$.

We can use Lemma~\ref{Lem:IntermediateRayPairSpecialCase} for this purpose: we have $S'=aS_k+\rho_\nu(r)$, then $S''=aS_k+\rho_\nu(\rho_\nu(r))$ etc.\ until the entries reach or exceed $aS_k+S_k$: if the entries reach $aS_k+S_k$, then $q_k=a+1$; if not, then $q_k>a+1$, and the only choice is $q_k=a+2$.
\qed

\proofof{Theorem~\ref{Thm:NumeratorsArbitrary} (Numerators arbitrary)}
We prove a stronger statement: given a hyperbolic component $W$ of period $n$, then we can combinatorially determine for any $s'>1$ how many components of period up to $s'$ are contained in the wake $\wake{W}$, how the wakes of all these components are nested, and what the widths of their wakes are; all we need to know about $W$ are the width of its wake and its internal address (without angles). In particular, these data encode how the internal address of the component $W$ can be continued within $\M$. This proves the theorem (and it also shows that the width of the wake of $W$ is determined by the internal address of $W$). 

For any period $s$, the number of periodic angles of period $s$ or dividing $s$ within any interval of $\Circle$ of length $\delta$ is either  $\lfloor\delta/(2^s-1)\rfloor$ or $\lceil\delta/(2^s-1)\rceil$ (the closest integers above and below $\delta/(2^s-1)$). If the interval of length $\delta$ is the wake of a hyperbolic component, then the corresponding parameter rays land in pairs, and the correct number of angles is the unique even integer among $\lfloor\delta/(2^s-1)\rfloor$ and $\lceil\delta/(2^s-1)\rceil$. This argument uniquely determines the exact number of hyperbolic components of any period within any wake of given width. 

There is a unique component $W_s$ of lowest period $s$, say, in $\wake{W}$ (if there were two such components, then this would imply the existence of at least one parameter ray of period less than $s$ and thus of a component of period less than $s$). The width of $\wake{W_{s}}$ is exactly $1/(2^s-1)$ (this is the minimal possible width, and greater widths would imply the existence of a component with period less than $s$).

Now suppose we know the number of components of periods up to $s'$ within the wake $\wake{W}$, together with the widths of all their wakes and how these wakes are nested. The wake boundaries of period up to $s'$ decompose $\wake{W}$ into finitely many components. Some of these components are wakes; the others are complements of wakes within other wakes. We can uniquely determine the number of components of period $s'+1$ within each of these wakes (using the widths of these wakes), and then also within each complementary component outside of some of the wakes (simply by calculating differences). Each wake and each complementary component can contain at most one component of period $s'+1$ by Theorem~\ref{Thm:CompletenessAngledIntAddr}. Note that the kneading sequences and hence the internal addresses of all wakes of period $s'+1$ are uniquely determined by those of periods up to $s'$. The long internal addresses tell us which wakes of period $s'+1$ contain which other wakes, and from this we can determine the widths of the wakes of period $s'+1$.
This provides all information for period $s'+1$, and this way we can continue inductively and prove the claimed statement.

Starting with the unique component of period $1$, it follows that the width of a wake $\wake{W}$ is determined uniquely by the internal address of $W$.
\qed

\proofof{Lemma~\ref{Lem:Numerators} (Finding numerators)}
\looseness -1
We only need to find the numerator $p_k$ if $q_k\ge 3$. Let $W_k$ be the unique hyperbolic component with angled internal address $(S_0)_{p_0/q_0}\IntAddr\dots\IntAddr(S_k)$ and let $\RP(\theta_k,\theta'_k)$ be the ray pair bounding its wake $\wake{W_k}$. Let $W'_k$ be the component with angled internal address $(S_0)_{p_0/q_0}\IntAddr\dots\IntAddr(S_k)_{p_k/q_k}\IntAddr q_kS_k$; it is an immediate bifurcation from $W_k$.

By Theorem~\ref{Thm:PropertiesParaRayPairs}, every $c\in \wake{W_k}$ has a repelling periodic point $z_c$ which is the landing point of the characteristic dynamic ray pair $\RP_c(\theta_k,\theta'_k)$; we find it convenient to describe our proof in such a dynamic plane, even though the result is purely combinatorial. Let $\Theta$ be the set of angles of rays landing at $z_c$; this is the same for all $c\in\wake{W_k}$. Especially if $c\in W'_k$, it is well known and easy to see that $\Theta$ contains exactly $q_k$ elements, the first return map of $z_c$ permutes the corresponding rays transitively and their combinatorial rotation number is $p_k/q_k$\version{ (Lemma~\ref{Lem:Permutation}) --- be more specific here!}{}. These rays disconnect $\C$ into $q_k$ sectors which can be labelled $V_0,V_1,\dots,V_{q_k-1}$ so that the first return map of $z_c$ sends $V_j$ homeomorphically onto $V_{j+1}$ for $j=1,2,\dots,q_{k}-2$, and so that $V_1$ contains the critical value and $V_0$ contains the critical point and the ray $R_c(0)$. Finally, under the first return map of $z_c$, points in $V_0$ near $z_c$ map into $V_1$, and points in $V_{q_k-1}$ near $z_c$ map into $V_0$. The number of sectors between $V_0$ and $V_1$ in the counterclockwise cyclic order at $z_c$ is then $p_k-1$, where $p_k$ is the numerator in the combinatorial rotation number. 
\hide{
These statements can also be proved directly for all $c\in\wake{W_k}$ using the defining property of characteristic ray pairs. 
}

Now suppose the dynamic ray $R_c(\theta)$ contains the critical value or lands at it. Then $R_c(\theta)\in V_1$, and $R_c(2^{(j-1)S_k}\theta)\in V_{j}$ for $j=1,2,\dots,q_k-1$. Counting the sectors between $V_0$ and $V_1$ means counting the rays $R_c(\theta),R_c(2^{S_k}\theta),\dots,R_c(2^{(q_k-2)S_k}\theta)$ in these sectors, and this means counting the angles $\theta, 2^{S_k}\theta, \dots, 2^{(q_k-2)S_k}$ in $(0,\theta)$. The numerator $p_k$ exceeds this number by one, and this equals the number of angles $\theta, 2^{S_k}\theta, \dots, 2^{(q_k-2)S_k}$ in $(0,\theta]$. 
\qed

\proofof{Lemma~\ref{Lem:LeftRightRay} (Left or right ray)}
The $n$-th entry in the kneading sequence of $\nu(\theta)$  is determined by the position of the angle $2^{n-1}\theta\in\{\theta/2,(\theta+1)/2\}$. The $n$-th entry in the kneading sequence of $W$ equals the $n$-th entry in the kneading sequence of $\tilde\theta$ for $\tilde\theta$ slightly greater than $\theta$ (for $\theta'$, we use an angle $\tilde\theta'$ slightly smaller than $\tilde\theta$); this is $\1$ if and only if $2^{n-1}\theta=\theta/2$ and $\0$ if $2^{n-1}\theta=(\theta+1)/2$. But $2^{n-1}\theta=\theta/2$ implies $2^{n-1}\theta\in (0,1/2)$, hence $b=0$, while $2^{n-1}\theta=(\theta+1)/2$ implies $2^{n-1}\theta\in(1/2,1)$ and $b=1$. The reasoning for $\theta'$ is similar.
\qed

\proofof{Proposition~\ref{Prop:IntAddrRenorm} (Internal Address and Renormalization)}
We only discuss the case of simple renormalization (the case of crossed renormalization is treated in \cite[Corollary~4.2]{CrossRenorm}).

Fix a hyperbolic component $W$ of period $n$. Let $\M_W$ be the component of $n$-renormalizable parameters in $\M$ containing $W$ (all $c\in W$ are $n$-renormalizable), and let $\Psi_W\colon\M\to\M_W$ be the tuning homeomorphism; see \cite{HaissinskyTuning,MiRenorm,MiOrbits} or \version{Section~\ref{Sec:Renormalization}}{\cite[\SecRenorm]{SymDyn}}. 
Let $1\IntAddr S_1\IntAddr\dots \IntAddr S_k=n$ be the internal address of $W$. Then the internal address of every $c\in\M_W$ starts with $1\IntAddr S_1\IntAddr\dots \IntAddr S_k=n$ because $\M_W$ contains no hyperbolic component of period less than $n$, so points in $\M_W$ are not separated from each other by parameter ray pairs of period less than $n$. All hyperbolic components within $\M_W$, and thus all ray pairs separating points in $\M_W$, have periods that are multiples of $n$, so all internal addresses of parameters within $\M_W$ have the form $1\IntAddr S_1\IntAddr\dots \IntAddr S_k\IntAddr S_{k+1}\dots$ so that all $S_m\ge n$ are divisible by $n$. In fact, if $c\in\M$ has internal address $1\IntAddr S'_1\IntAddr\dots \IntAddr S'_{k'}\dots$, then it is not hard to see that the internal address of $\Psi_W(c)$ is $1\IntAddr S_1\IntAddr\dots \IntAddr n\IntAddr nS'_1\IntAddr\dots \IntAddr nS'_k\dots$ (hyperbolic components of period $S'$ in $\M$ map to hyperbolic components of period $nS'$ in $\M_W$, and all ray pairs separating points in $\M_W$ are associated to hyperbolic components in $\M_W$ that are images under $\Psi_W$).

For the converse, we need \emph{dyadic Misiurewicz-Thurston parameters}: these are by definition the landing points of parameter rays $R(\theta)$ with $\theta=m/2^k$; dynamically, these are the parameters for which the singular orbit is strictly preperiodic and terminates at the $\beta$ fixed point. If $\theta=m/2^k$, then the kneading sequence $\nu(\theta)$ has only entries $\0$ from position $k+1$, so the internal address of $\nu(\theta)$ contains all integers that are at least $2k-1$. But the internal address of $\nu(\theta)$ from Algorithm~\ref{Def:IntAddrKneading} equals the internal address of $R(\theta)$ in parameter space (Proposition~\ref{Prop:EqualIntAddr}), and the landing point of $R(\theta)$ has the same internal address. Therefore, the internal address of any dyadic Misiurewicz-Thurston parameter  contains all sufficiently large positive integers.

Suppose the internal address of some $c\in\M$ has the form $1\IntAddr S_1\IntAddr\dots \IntAddr S_k\IntAddr S_{k+1}\dots$ with $S_k=n$ and all $S_m\ge n$ are divisible by $n$. There is a component $W$ with address $1\IntAddr S_1\IntAddr\dots \IntAddr S_k$ so that $c\in\wake{W}$.
If $c\not\in\M_W$, then $c$ is separated from $\M_W$ by a Misiurewicz-Thurston parameter $c_*\in\M_W$ which is the tuning image of a dyadic Misiurewicz-Thurston parameter (see \cite{DoAngles}, \cite[Section~8]{MiOrbits}, or \version{Corollary~\ref{Cor:BoundaryRenormalization}}{\cite[\CorBoundaryRenormalization]{SymDyn}}). Therefore, the internal address of $c_*$ contains all  integers that are divisible by $n$ and sufficiently large, say at least $Kn$. 
Let $S$ be the first entry in the internal address of $c$ that corresponds to a component ``behind $\M_W$'' (so that it is separated from $\M_W$ by $c_*$). The long internal address of $c$ contains ``behind'' $W$ only hyperbolic components of periods divisible by $n$, and this implies that before and after $c_*$ there must be two components of equal period (greater than $Kn$) that are not separated by a ray pair of lower period. This contradicts Lavaurs' Lemma\version{~\ref{Lem:Lavaurs}}{}.
\qed

\remark
In Definition~\ref{Def:AngledIntAddr}, we defined angled internal addresses in parameter space, but they can also be defined dynamically: the angles are combinatorial rotation numbers of rays landing at characteristic periodic points the periods of which occur in the internal address. This allows us to give more dynamic proofs of several theorems that we proved in parameter space. For instance, changing numerators has no impact on whether an angled internal address is realized in the complex plane (Theorem~\ref{Thm:NumeratorsArbitrary}): for the Hubbard trees at centers of hyperbolic components, this simply changes the embedding of the tree, but not the issue whether such a tree can be embedded (compare the discussion in \version{Corollary~\ref{CorNumberEmbed} and Lemma~\ref{LemNumberEmbed2})}{\cite[{\CorNumberEmbed} and {\LemNumberEmbedTwo}]{SymDyn})}.  (The denominators are determined already by the internal address without angles; see Lemma~\ref{Lem:Denominators}). Similarly, the angled internal address completely determines a Hubbard tree with its embedding and thus the dynamics at the center of a hyperbolic component; this is Theorem~\ref{Thm:CompletenessAngledIntAddr}: the (finite) angled internal address completely specifies every hyperbolic component.

\version{
We can define a partial order of hyperbolic components $W,W'$ in $\M$ as follows: we say that $W'$ is greater than $W$ if $W'\neq W$ and the wake of $W$ contains $W'$ (and thus the wake of $W'$). Since wakes are by definition bounded by parameter ray pairs at periodic angles, and are thus either nested or disjoint, Theorem~\ref{Thm:PropertiesParaRayPairs} implies that an equivalent dynamical definition goes as follows: $W'$ is greater than $W$ if any parameter $c'\in W'$ (or even $c'$ in the wake of $W'$) has the property that if $\RP(\theta,\theta')$ is the ray pair bounding the wake of $W$, then $p_{c'}$ has a characteristic dynamic ray pair $\RP_{c'}(\theta,\theta')$.
(An extension of this partial order to all $c,c'\in\M$ goes as follows: $c'$ is greater than $c$ if there is a parameter ray pair $\RP(\theta,\theta')$ at periodic angles that separates $c'$ from $c$ and the origin.)
}{}

\version{
In Definition~\ref{DefOrder}, the space of all $\*$-periodic kneading sequences (realized by a complex polynomial or not) was endowed with a partial order, so that $\nu'>\nu$ if the Hubbard tree for $\nu'$ contains a periodic point with itinerary $\A(\nu)$. These two partial orders are related as follows.
}{}

\version{
\begin{prop}{Order for Hyperbolic Components and Kneadings}
\label{Prop:OrderCompsKneadings} \lineclear
Consider two hyperbolic components $W\neq W'$ with associated kneading sequences $\nu$ and $\nu'$. If $W'$ is greater than $W$, then $\nu'>\nu$. Conversely, if $\nu'>\nu$, then there is a hyperbolic component $W''$ greater than $W$ so that $W''$ has kneading sequence $\nu'$ and thus the same period and internal address as $W'$.
\qedd
\end{prop}
We omit the proof, which is not difficult: it follows from the dynamical definition of the partial order as given above.
Note that it is not true that if $\nu'>\nu$, then $W'$ is greater than $W$: there can be several hyperbolic components with the same kneading sequence $\nu'$ and thus the same internal address, and these components are distinguished by their angled internal addresses. At least one of these components is greater than $W$, while the other components may not be comparable.
}{}


\version
{\newpage\section{Symbolic Dynamics, Permutations, Galois Groups}}
{\section{Narrow Components and Admissibility}}

\label{Sec:Permutations}

Angled internal addresses describe hyperbolic components, or more generally combinatorial classes, in $\M$ uniquely. But are all possible angled internal addresses indeed realized (``admissible'')? Since the angles have no impact on admissibility (Theorem~\ref{Thm:NumeratorsArbitrary}), this depends only on the internal address without angles. An infinite internal address is admissible if and only if all its finite truncations are admissible (because $\M$ is closed), so it suffices to investigate admissibility of all finite internal addresses $1\IntAddr \dots \IntAddr S_k$. An equivalent question is whether all sequences $\ovl{\1 \nu_2\nu_3\dots \nu_{n-1}\*}$ with $\nu_i\in\{\0,\1\}$ are kneading sequences of periodic external angles. 

Milnor and Thurston \cite{MilnorThurston} classified all \emph{real-admissible} kneading sequences, i.e., those sequences that occur for real quadratic polynomials: an $n$-periodic sequence in $\{\0,\1\}^\Nzero$ is real-admissible if and only if it is the largest in a certain (non-lexicographic) order among its $n$ periodic shifts; so roughly a fraction of $1/n$ of all $n$-periodic sequences is real-admissible (note the analogy to our Lemma~\ref{Lem:MaxShift}: if any periodic sequence is the largest, among all its shifts, with respect to lexicographic order, then it corresponds to a certain kind of (usually non-real) hyperbolic components that we call ``purely narrow'').

It turns out that not all sequences are (complex-)admissible either. This can be seen already by a statistical argument: there are $2^{n-2}$ possible kneading sequences of period $n$ (starting the period with $\1$ and ending with $\*$), and there are $2^n$ external rays of period $n$, landing in pairs at $2^{n-1}$ hyperbolic components of period $n$ (both up to terms of order $O(2^{n/2})$ corresponding to exact periods strictly dividing $n$). So every possible kneading sequence of period $n$ should be realized \emph{on average} by two hyperbolic components. Those on the real axis are only a small fraction, and all the others appear in pairs corresponding to complex conjugation of $\M$ (or equivalently replacing in the angled internal address the first angle $p_k/q_k\neq 1/2$ with $1-p_k/q_k$). But a significant fraction of internal addresses exists more than twice: whenever there is a denominator $q_k\not\in\{2,3,4,6\}$, then there are more than two numerators possible, or also when there are at least two denominators other than $1/2$. The average time an admissible internal address is realized is thus greater than $2$, so that the non complex-admissible internal addresses have positive measure (among all $\0$-$\1$-sequences with the product measure). In \cite{AdmissKneading}, a complete classification of all complex admissible kneading sequences was given, but the criterion is rather non-trivial; and in \cite{BS-Admiss-PosMeasure} it was shown that complex admissible kneading sequences have positive measure among all $\0$-$\1$-sequences.

The simplest candidate internal address that is not complex admissible is $1\IntAddr 2\IntAddr 4\IntAddr 5\IntAddr 6$ with kneading sequence $\nu=\ovl{\1\0\1\1\0\0}$. If it existed, the kneading sequence ``just before'' the component would be $\Abar(\nu)=\ovl{\1\0\1\1\0\1}$ of period $3$, so the period $6$ component (the limiting kneading sequence outside the wake) would be a bifurcation from a period $3$ component, but the period $3$ component does not show up in the internal address. Moreover, the period $3$ component would have kneading sequence $\ovl{\1\0\1}$ with an infinite internal address! There is a hyperbolic component of period $3$ with internal address $\1\IntAddr 2\IntAddr 3$ and kneading sequence $\ovl{\1\0\0}$, so our candidate period $6$ component would bifurcate from this component ``outside of the wake''. So this internal address does not exist, and neither does any that starts with $1\IntAddr 2\IntAddr 4\IntAddr 5\IntAddr 6\IntAddr\dots$
It turns out \cite{AdmissKneading} that any non-existing kneading sequence is generated by a similar process.

How can one describe all admissible internal addresses? Given any finite internal address $1\IntAddr S_1\IntAddr \dots \IntAddr S_k$, what are the possible continuations $S_{k+1}>S_k$? In this section, we give a precise answer for a particular kind of components that we call ``purely narrow": a ``narrow'' component is one for which we have good control, and purely narrow means that we construct the internal address along components that are all narrow.

\version{ 
\Intro{We investigate which permutations of periodic points of $p_c$ can be achieved by analytic continuation with respect to the parameter $c$. This makes it possible to determine Galois groups of those polynomials which represent the periodic points. In order to do this, we need some existence theorems of kneading sequences, and these will be derived from results on the structure of the Mandelbrot set using internal addresses.} }{}

Every component $W$ of period $n$ has an associated internal address (finite, ending with entry $n$); compare the remark before Definition~\ref{Def:AngledIntAddr}. The component $W$ also has an associated periodic kneading sequence $\nu(W)$ consisting only of entries $\0$ and $\1$: one way of defining this is to take any parameter ray $R(\phi)$ with irrational $\ph$ landing at $\partial W$; then $\nu(W):=\nu(\phi)$. Equivalently, let $\RP(\theta,\theta')$ be the parameter ray pair landing at the root of $W$; then $\theta$ and $\theta'$ have period $n$ and $\nu(\theta)=\nu(\theta')$ are $\*$-periodic of period $n$, and $\nu(W)=\A(\nu(\theta))=\A(\nu(\theta'))$ (see Lemma~\ref{Lem:ParaRayKneading} and Proposition~\ref{Prop:EqualIntAddr}). Of course, the internal address of $W$ is the same as the internal address of $R(\phi)$, $R(\theta)$, $R(\theta')$ or of $\phi$, $\theta$ or $\theta'$.

\begin{defi}{Narrow Component}
\label{Def:Narrow} \lineclear
A hyperbolic component of period $n$ is \emph{narrow} if its wake contains no component of lower period, or equivalently if the wake has width $1/(2^n-1)$.
\end{defi}

\remark
It follows directly from 
\version{Proposition~\ref{Prop:WidthWake}}{(\ref{Eq:WidthWake})} that if $W'$ is a bifurcation from another component $W$, then $W'$ is narrow if and only if $W$ has period $1$. 

If $W,W'$ are two hyperbolic components of periods $n$ and $n'$ so that $W'$ is contained in the wake $\wake{W}$ of $W$, then we say that $W'$ is \emph{visible} from $W$ if there exists no parameter ray pair of period less than $n'$ that separates $W$ and $W'$. If $n<n'$, then this is equivalent to the condition that the internal address of $W'$ be formed by the internal address of $W$, directly extended by the entry $n'$.

\begin{lemma}{Visible Components from Narrow Component}
\label{Lem:VisibleNarrow} \lineclear
For every narrow hyperbolic component of period $n$, there are visible components of all periods greater than $n$. More precisely, every $p/q$ sublimb contains exactly $n$ visible components: exactly one component each of period $qn-(n-1)$, $qn-(n-2)$,\dots, $qn$.
\end{lemma}
\proof 
Let $W$ be a narrow hyperbolic component of period $n$ and consider its $p/q$-subwake.
The visible components in this wake have periods at most $qn$. First we show that any two
visible components within the $p/q$-subwake of $W$ have different internal addresses.
By way of contradiction, suppose there are two
components $W_1$ and $W_2$ of equal period $m\le qn$ with the same internal 
address. By Theorem~\ref{Thm:CompletenessAngledIntAddr} there must
be another hyperbolic component $W'$ in the same subwake of $W$ so
that $W_1$ and $W_2$ are in different $p_1/q'$- and $p_2/q'$-subwakes of $W'$
with $q'\ge 3$ (different hyperbolic components with the same internal address
must have different angled internal addresses). Let $n'> n$ be the period of $W'$. 
Then the width of these subwakes is, according to 
\version{Proposition~\ref{Prop:WidthWake},}{(\ref{Eq:WidthWake}),}
\begin{eqnarray*}
|W'|\frac{(2^{n'}-1)^2}{2^{q'n'}-1}
&\le&
 |W| \frac{(2^n-1)^2}{2^{qn}-1}\cdot\frac{(2^{n'}-1)^2}{2^{q'n'}-1}
= \frac{2^n-1}{2^{qn}-1}\cdot\frac{(2^{n'}-1)^2}{2^{q'n'}-1} 
\\
&<& \frac{2^n(2^{(2-q')n'})(1-2^{-q'n'})^{-1}} {2^{qn}-1}
\le \frac{2^{(n-n')} \cdot 2}  {2^{qn}-1}
\le \frac{1}{2^{qn}-1}
\end{eqnarray*}
because $W'$ is contained in the $p/q$-wake of $W$. But this is not large enough to contain a component of period at most $qn$, a contradiction. 

The next step is to prove that every subwake of $W$ of denominator $q$ contains
exactly $2^k$ external rays with angles $a/(2^{(q-1)n+k}-1)$ for
$1\leq k\leq n$, including the two rays bounding the wake. In fact,
\version{Proposition~\ref{Prop:WidthWake} says that }{} the width of the wake is
$(2^n-1)/(2^{qn}-1)$, so the number of rays one expects by comparing
widths of wakes is
\[
\frac{(2^n-1)(2^{(q-1)n+k}-1)}{2^{qn}-1}= 2^k+\frac{2^k-2^n-2^{(q-1)n+k}+1}{2^{qn}-1} =: 2^k+\alpha,
\] 
where an easy calculation shows that $-1\leq\alpha<1$ and that
$\alpha=-1$ can occur only for $k=n$. The actual number of rays can
differ from this expected value by no more than one and is even, hence
equal to $2^k$. Moreover, no such ray of angle $a/(2^{(q-1)n+k}-1)$
can have period smaller than $(q-1)n+k$ because otherwise it would land at a
hyperbolic component of some period dividing $(q-1)n+k$; but in the
considered wake  there would not be room enough to contain a second
ray of equal period.

This shows that, for any $k\leq n$, the number of hyperbolic
components of period $m=(q-1)n+k$ in any subwake of $W$ of
denominator $q$  equals $2^{k-1}$. 
They must all have different internal addresses.
The single component of period $(q-1)n+1$ takes care of the case $k=1$, and its wake subdivides the $p/q$-subwake of $W$ into two components. There are two components of period $(q-1)n+2$, and since their internal addresses are different, exactly one of them must be in the wake of the component of period $(q-1)n+1$, while the other is not; the latter one is visible from $W$. (The non-visible component is necessarily narrow, while the visible component may or may not be narrow.)

So far we have taken are of $3$ components, and their wake boundaries subdivide the $p/q$-subwake of $W$ into $4$ pieces. Each piece most contain one component of period $(q-1)n+3$, so exactly one of these components is visible from $W$, and so on. The argument continues as long as we have uniqueness of components for given internal addresses, which is for $k\le n$.
\qed


\begin{lemma}{Narrow Visible Components from Narrow Component}
\label{Lem:NarrowVisibleNarrow} \lineclear
Suppose $W$ and $W'$ are two hyperbolic components of periods $n$ and $n+s$ with $s>0$ so that $W$ is narrow and $W'$ is visible from $W$. Let $k\in\{1,\dots,n-1,n\}$ be so that $s\equiv k$ modulo $n$. 
Then the question whether or not $W'$ is narrow depends only on the first $k$ entries in the kneading sequence of $W$ (but not otherwise on $W$).
\end{lemma}
\proof
By Lemma~\ref{Lem:VisibleNarrow}, every $p/q$-subwake of $W$ contains exactly one visible hyperbolic component $W_m$ of period $m=(q-1)n+k$ for $k=1,2,\dots,n$. Such a component $W_m$ is narrow unless the wake $\wake{W_m}$ contains one of the components $W_{m'}$ with $m'=(q-1)n+k'$ and $k'\in\{1,2,\dots,k-1\}$.
If there is such a component, let $m'$ be so that the width $|W_{m'}|$ is maximal; then $W_m$ and $W_{m'}$ cannot be separated by a parameter ray pair of period less than $m'$. This makes $m'$ unique.

In order to find out whether the visible component $W_{m'}$ is contained in $\wake{W_m}$, we need to compare the kneading sequence $\nu$ associated to $W$ with the kneading sequence $\nu'$ ``just before'' $W_{m'}$ and determine whether their first difference occurs at position $m$: according to Corollary~\ref{Cor:IntermediateRayPair}, the position of the first difference is the least period of two ray pairs separating $W$ and $W_{m'}$. If this difference occurs after entry $m$, then the ray pairs landing at the root of $W_m$ do not separate $W$ from $W_{m'}$; the difference cannot occur before entry $m$ because of visibility of $W_m$ and the choice of $m'$. 

By visibility of $W_{m'}$ from $W$, the kneading sequences $\nu$ and $\nu'$ coincide for at least $m'>n$ entries. Eliminating these, we need to compare $\sigma^{m'}(\nu)$ with $\sigma^{m'}(\nu')=\nu'$ for $m-m'=k-k'<n$ entries. 
But the first $k-k'$ entries in $\sigma^{m'}(\nu)$ coincide with those in $\sigma^{k'}(\nu)$ (because $\nu$ has period $n$), and the first $k-k'$ entries in $\nu'$ and $\nu$ are also the same; so we need to compare the first $k-k'$ entries in $\sigma^{k'}(\nu)$ with $\nu$: these are the entries $\nu_{k'+1}\dots\nu_{k}$ and $\nu_1\dots\nu_{k-k'}$. This comparison involves only the first $k$ entries in $\nu$. (The precise criterion is: $W'$ is not narrow if and only if there is a $k'\in\{1,2,\dots,k-1\}$ with $\rho(k')=k$%
\version{; compare (\ref{Eq:no_jump}) in the proof of Corollary~\ref{Cor:SufficientAdmiss}}{}.) \comment{Give ref? See note in source file}
\qed


\begin{prop}{Combinatorics of Purely Narrow Components}
\label{Prop:CombinatoricsPurelyNarrow} \lineclear
Consider a hyperbolic component $W$ with internal address $1\IntAddr S_1\IntAddr \dots \IntAddr S_k$ and associated kneading sequence $\nu$. Suppose that $\nu_{S_i}=\0$ for all $i=1=2,\dots,k$. Then $W$ is narrow and for every $S_{k+1}>S_k$, there exists a hyperbolic component with internal address $1\IntAddr S_1\IntAddr \dots \IntAddr S_k\IntAddr S_{k+1}$; it is narrow if and only if $\nu_{S_{k+1}}=\1$.
\end{prop}

\proof
We prove the claim by induction on the length of internal addresses, starting with the address $1$ of length $0$. The associated component has period $1$, it is narrow and has $\nu=\ovl{\1}$, and there is no condition on $\nu_{S_{i}}$ to check. For every $S_1>1$, there exists a hyperbolic component with internal address $1\IntAddr S_1$ (these are components of period $S_1$ bifurcating immediately from the main cardioid). By the remark after Definition~\ref{Def:Narrow}, these components are narrow, and indeed $\nu_{S_1}=\1$.

Now assume by induction that the claim is true for a narrow component $W_{k-1}$ with internal addresses $1\IntAddr S_1\IntAddr \dots \IntAddr S_{k-1}$ of length $k-1$ and associated kneading sequence $\mu$. We will prove the claim when $W$ is a hyperbolic component with internal address $1\IntAddr S_1\IntAddr \dots \IntAddr S_{k-1}\IntAddr S_k$ of length $k$ and with associated internal address $\nu$ so that $\nu_{S_k}=\0$. First we show that $W$ is narrow: by the inductive hypothesis, it is narrow if and only if $\mu_{S_k}=\1$, but this is equivalent to $\nu_{S_k}=\0$ because $\nu$ and $\mu$ first differ at position $S_k$.

Consider some integer $S_{k+1}>S_k$. 
By Lemma~\ref{Lem:VisibleNarrow}, there exists another component $W_{k+1}$ with internal address $1\IntAddr S_1\IntAddr \dots \IntAddr S_k\IntAddr S_{k+1}$. We need to show that it is narrow if and only if  $\nu_{S_{k+1}}=\1$. 
If $S_{k+1}$ is a proper multiple of $S_k$, then the assumption $\nu_{S_k}=\0$ implies $\nu_{S_{k+1}}=\0$, and we have to show that $W_{k+1}$ is not narrow. Indeed, $W_{k+1}$ is a bifurcation from $W$ by Lemma~\ref{Lem:EstimatesDenominators}, and by the remark after Definition~\ref{Def:Narrow} bifurcations are narrow if and only if they bifurcate from the period $1$ component. 

We can thus write $S_{k+1}=qS_k+S'_{k+1}$ with $S'_{k+1}\in\{S_k+1,S_k+2,\dots,2S_k-1\}$ and $q\in\Nzero$. 
Again by Lemma~\ref{Lem:VisibleNarrow}, there exists another component $W'_{k+1}$ with internal address $1\IntAddr S_1\IntAddr \dots \IntAddr S_k\IntAddr S'_{k+1}$; by Lemma~\ref{Lem:NarrowVisibleNarrow}, it is narrow if and only if $W_{k+1}$ is. 
But $\nu$ has period $S_k$, so $\nu_{S_{k+1}}=\nu_{S'_{k+1}}$ and the claim holds for $S_{k+1}$ if and only if it holds for $S'_{k+1}$. It thus suffices to restrict attention to the case $S_{k+1}<2S_k$.
Whether or not $W_{k+1}$ is narrow is determined by the initial $S_{k+1}-S_k<S_k$ entries in $\nu$.

Now we use the inductive hypothesis for $W_{k-1}$: there exists a component $W'$ with address $1\IntAddr S_1\IntAddr \dots \IntAddr S_{k-1}\IntAddr (S_{k-1}+S_{k+1}-S_k)$, and it is narrow if and only if $\mu_{S_{k-1}+S_{k+1}-S_k}=\1$. The kneading sequences $\mu$ and $\nu$ first differ at position $S_k$, so their initial $S_{k+1}-S_k$ entries coincide. Therefore, again by Lemma~\ref{Lem:NarrowVisibleNarrow}, the component $W'$ is narrow if and only if $W_{k+1}$ is. Therefore, $W_{k+1}$ is narrow if and only if $\mu_{S_{k-1}+S_{k+1}-S_k}=\1$.  Finally, we have
\[
\nu_{S_{k+1}}=\nu_{S_{k+1}-S_k}=\mu_{S_{k+1}-S_k}=\mu_{S_{k-1}+S_{k+1}-S_k}
\]
by periodicity of $\nu$ (period $S_k$) and of $\mu$ (period $S_{k-1}$) and because the first difference between $\nu$ and $\mu$ occurs at position $S_k>S_{k+1}-S_k$. This proves the proposition.
\qed

\remark
We call a hyperbolic component $W_{k+1}$ with internal address $1\IntAddr S_1\IntAddr \dots \IntAddr S_k\IntAddr S_{k+1}$ and associated kneading sequence $\nu$ \emph{purely narrow} if $\nu_{S_i}=\0$ for $i=1,2,\dots,k+1$. Proposition~\ref{Prop:CombinatoricsPurelyNarrow} implies that this is equivalent to the condition that $1\IntAddr S_1\IntAddr \dots \IntAddr S_{i-1}\IntAddr S_{i}$ describes a narrow component for $i=1,\dots,k+1$ (hence the name). The asymmetry in the statement of the proposition (for narrow components, the last entry in $\nu$ must be $\1$, rather than $\0$ for all earlier components) is because the condition is with respect to the kneading sequence of period $S_k$, not with respect to the sequence of period $S_{k+1}$ associated to $W_{k+1}$.

\begin{figure}[ptb]
\includegraphics{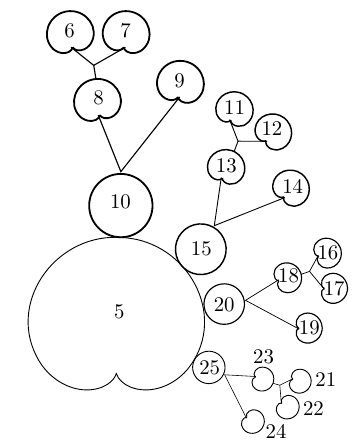}
\includegraphics[trim=5 0 0 0]{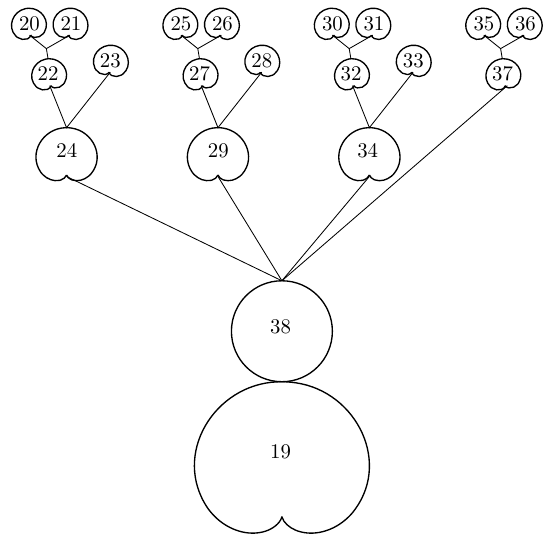}
\\ \begin{picture}(0,0)
\put(35,43){\vector(-4,1){90}}
\end{picture}
\caption{Left: For a narrow hyperbolic component $W$ (here of period $5$), the trees of visible components within any $p/q$- and $p'/q'$-subwake are the same when adding $(q'-q)n$ to the periods of all components in the $p/q$-subwake; this ``translation principle'' follows directly from Lemma~\ref{Lem:NarrowVisibleNarrow}, using just the combinatorics of internal addresses or kneading sequences (even the embeddings of the trees are the same;  this follows from comparisons with dynamical planes). Right: if $W'$ (here of period $19$) is visible from $W$ and both are narrow, then the tree of visible components within the $1/2$-subwake of $W'$ can be reconstructed from the trees of visible components within various subwakes of $W$: if $n'$ and $n$ are the periods of $W'$ and $W$, then the tree formed by the visible components in the $1/2$-subwake of $W'$ of periods $n'+1$, \dots, $2n'-1$ (excluding the bifurcating component of period $n'$) equals the tree formed by the visible components of periods $n+1,\dots,n+n'-1$ in the $1/q$-subwakes of $W$, adding $n'-n$ to all periods.}
\label{Fig:Translation&Growing}
\end{figure}
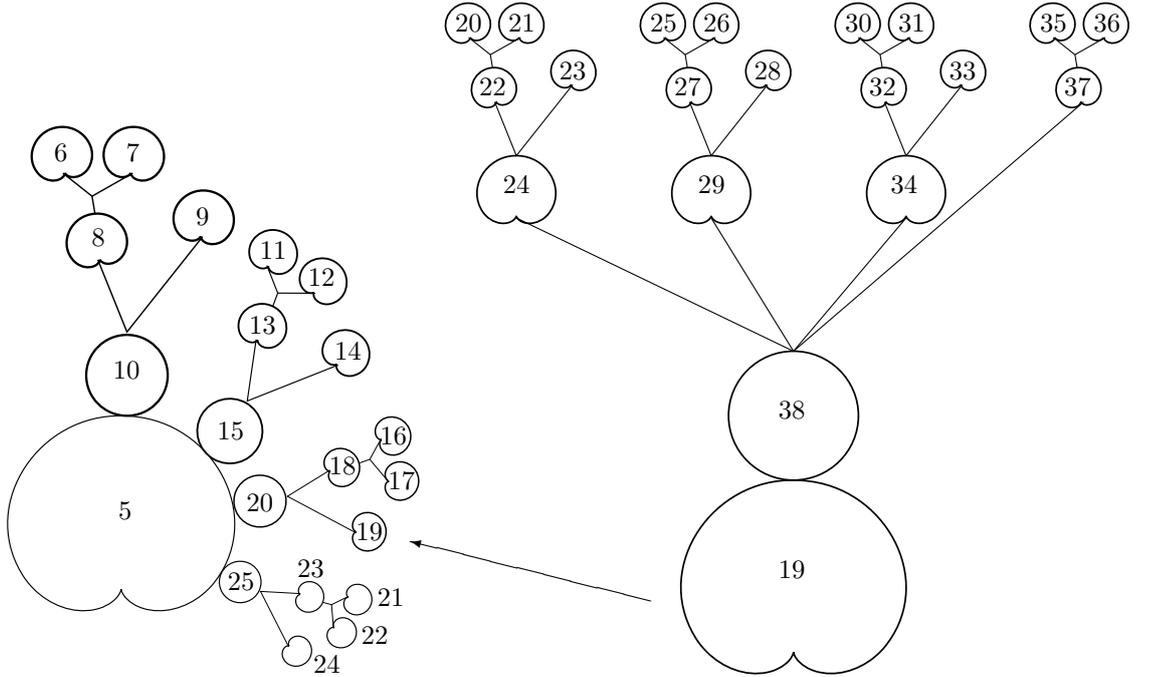

\remark
For every narrow hyperbolic component, the previous results allow to construct combinatorially the trees of visible components within any sublimb; for purely narrow components, the global tree structure can thus be reconstructed by what we call ``growing of trees'': see Figure~\ref{Fig:Translation&Growing}. 
These issues have been explored further by Kauko~\cite{VirpiPaper}.

\begin{lemma}{Maximal Shift of Kneading Sequence}
\label{Lem:MaxShift} \lineclear
For a kneading sequence $\nu$ (without $\*$) with associated internal address $1\IntAddr S_1\IntAddr \dots \IntAddr S_k\IntAddr\dots$ (finite or infinite), the following are equivalent:
\begin{enumerate}
\item
no shift $\sigma^k(\nu)$ exceeds $\nu$ with respect to lexicographic ordering;
\item
for every $r\ge 1$, we have $\nu_{\rho(r)}=\0$.
\item
for every $k\ge 1$, we have $\nu_{S_k}=\0$.
\end{enumerate}
In all these cases, $r<S_k$ implies $\rho(r)\le S_k$.
\end{lemma}
\proof
Pick some $r\ge 1$. 
Comparing $\sigma^r(\nu)$ to $\nu$ in the lexicographic ordering means comparing the entries $\nu_{r+1}\nu_{r+2}\dots\nu_{\rho(r)}$ to the entries $\nu_1\nu_2\dots\nu_{\rho(r)-r}$, where the position of first difference occurs by definition at position $\rho(r)$ (in $\sigma^r(\nu)$) vs.\ at position $\rho(r)-r$ (in $\nu$). Therefore, $\nu>\sigma^r(\nu)$ if and only if $\nu_{\rho(r)}=\0$. Thus the first two conditions are  equivalent, and the second implies the third because each $S_k=\rho(S_{k-1})$.

Now suppose that $\nu$ satisfies the third condition. For $k=1,2,\dots$, define $\nu^k$ to be the kneading sequence (without $\star$) corresponding to the finite internal address $1\IntAddr S_1\IntAddr \dots\IntAddr S_k$, so that $\nu^k$ has period $S_k$. Let $\rho_{\nu^k}$ be the non-periodicity function associated to $\nu^k$, i.e., the first difference between $\sigma^r(\nu^k)$ and $\nu^k$ occurs at position $\rho_{\nu^k}(r)-r$.

We will show by induction that for all $r< S_k$, we have $\rho_{\nu^k}(r)\le S_k$ and ${\nu^{k}}_{\rho_{\nu^{k}}(r)}=\0$, assuming that this holds for $S_{k-1}$ and $\nu^{k-1}$. \hide{Write $S_k=qS_{k-1}+s$. }

Suppose there is an $r$ so that $\rho_{\nu^k}(r)>S_k$. 
We have $\rho_{\nu^{k-1}}(r)=S_k$ because $\nu^{k-1}$ and $\nu^k$ differ at position $S_k$ (for the first time). 
By hypothesis, ${\nu^k}_{S_k}=\nu_{S_k}=\0$, hence ${\nu^{k-1}}_{S_k}=\1$.
But now we have ${\nu^{k-1}}_{\rho_{\nu^{k-1}}(r)}=\1$, in contradiction to the inductive hypothesis.
\hide{
hence ${\nu^{k-1}}_{\rho(r'')}=\1$ for $r''=r-q'S_{k-1}$ in contradiction to the inductive hypothesis.
}

Now we show ${\nu^k}_{\rho_{\nu^k}(r)}=\0$ for all $r$. If $\rho_{\nu^k}(r)<S_k$, then we have $\rho_{\nu^k}(r)=\rho_{\nu^{k-1}}(r)$ and thus ${\nu^{k}}_{\rho_{\nu^k}(r)}={\nu^{k-1}}_{\rho_{\nu^{k-1}}(r)}=\0$ by inductive hypothesis. And if $\rho_{\nu^k}(r)=S_k$, then the claim holds by hypothesis.

Therefore condition (3) implies (2), so all three are equivalent.
\qed

\begin{prop}{Internal Address of Purely Narrow Hyp.\ Component}
\lineclear
Every finite internal address of the type described in Lemma~\ref{Lem:MaxShift} is realized and corresponds to purely narrow hyperbolic components, and conversely each purely narrow hyperbolic component has an internal address of this type. 
\end{prop}
\proof
A hyperbolic component with internal address $1\IntAddr S_2\IntAddr \dots \IntAddr S_k$ and associated kneading sequence $\nu$ is purely narrow by definition if $\nu_{S_k}=\0$ for all $k$.

Conversely, it follows inductively from Proposition~\ref{Prop:CombinatoricsPurelyNarrow} that all internal addresses $1\IntAddr S_1\IntAddr \dots \IntAddr S_k$ are realized by narrow components because the associated kneading sequence $\nu$ satisfies $\nu_{S_i}=\nu_{\rho(S_{i-1})}=\0$.
\qed

\hide{
\begin{lemma}{Maximal Shift of Kneading Sequence}
\label{Lem:MaxShift} \lineclear
For a kneading sequence $\nu$ (without $\*$) with associated internal address $1\IntAddr S_1\IntAddr \dots \IntAddr S_k\IntAddr\dots$ (finite or infinite), the following are equivalent:
\begin{enumerate}
\item
no shift $\sigma^k(\nu)$ exceeds $\nu$ with respect to lexicographic ordering;
\item
for every $r\ge 1$, $\nu_{\rho(r)}=\0$.
\end{enumerate}
Every finite internal address of this type is realized by a purely narrow hyperbolic component. \comment{Converse?}
\end{lemma}
\proof
Pick some $r\ge 1$. If $\nu_{\rho(r)}=\1$, then the entries $\nu_{r+1}\nu_{r+2}\dots\nu_{\rho(r)}$ exceed the entries $\nu_1\nu_2\dots\nu_{\rho(r)-r}$, hence $\sigma^{r}(\nu)>\nu$ in the lexicographic ordering; conversely, if $\nu_{\rho(r)}=\0$, then $\sigma^{r}(\nu)<\nu$. Thus both conditions are indeed equivalent.
By Proposition~\ref{Prop:CombinatoricsPurelyNarrow}, it follows inductively that all internal addresses $1\IntAddr S_1\IntAddr \dots \IntAddr S_{k-1}\IntAddr S_k$ are realized by narrow components because the associated kneading sequence $\nu$ satisfies $\nu_{S_i}=\nu_{\rho(S_{i-1})}=\0$.
\qed
}

\remark
The existence of kneading sequences as described in this result can also be derived from the general admissibility condition on kneading sequences, see \version{Corollary~\ref{Cor:SufficientAdmiss}}{\cite[\CorMaxShiftAdmiss]{SymDyn}} (this is a more abstract and difficult, but also more general result).
Schmeling observed that sequences that are maximal shifts with respect to the lexicographic order are exactly the fixed points of the map $\theta\mapsto\nu(\theta)$. Recently, Buff and Tan Lei observed \cite{BuffTanLei} that this gives a simple combinatorial proof of the existence of kneading sequences as in Lemma~\ref{Lem:MaxShift}, without interpretation in parameter space.

\version
{\heading{Symbolic Dynamics and Permutations}}
{\subsection{Symbolic Dynamics and Permutations}}
\label{Sub:Permutations}

We will now discuss permutations of periodic points of $p_c(z)=z^2+c$. We make a brief excursion to algebra and describe our theorem first in algebraic terms (for readers that are less familiar with these algebraic formulations, we restate the result  in Theorem~\ref{Thm:PermutationPeriodicPoints}). For $n\ge 1$, let $Q_n(z):=p_c^{\circ n}(z)-z$ (consider these as polynomials in $z$ with coefficients in $\C[c]$). The roots of $Q_n$ are periodic points of period dividing $n$, so we can factor them as
\[
Q_n = \prod_{k|n} P_k \,\,;
\]
this product defines the $P_k$ recursively, starting with $P_1=Q_1$. 
\begin{theo}{Galois Groups of Polynomials}
\label{Thm:GaloisGroups} \lineclear
For every $n\ge 1$, the polynomials $P_n$ are irreducible over $\C[c]$. Their Galois groups $G_n$ consist of all the permutations of the roots of $P_n$ that commute with the dynamics of $p_c$. There is a short exact sequence
\[
0 \longrightarrow (\Z_n)^{N_n} \longrightarrow G_n \longrightarrow S_{N_n} \longrightarrow 0
\,\,,
\]
where $\Z_n=\Z/n\Z$, while $N_n$ is the number of periodic orbits of exact period $n$ for $p_c$ with $c\in X_n$, and $S_{N_n}$ is the symmetric group on $N_n$ elements.
\end{theo}
In this statement, the injections $(\Z_n)^{N_n}\to G_n$ correspond to independent cyclic permutations of the $N_n$ orbits of period $n$, while the surjection is the projection from periodic points to periodic orbits and yields arbitrary permutations among the orbits. This theorem was first shown by Bousch \cite{Bousch} by more algebraic methods and more recently by Morton and Patel in \cite{MP}.

A related statement in parameter space is still unsolved: consider the polynomials $\tilde Q_n(c):=Q_n(c)=p_c^{\circ n}(c)-c \in\Z[c]$. Their roots are parameters $c$ for which the critical orbit is periodic of period dividing $n$ (i.e., $c$ is the center of a hyperbolic component of period $n$), so we can again factor as $\tilde Q_n=\prod_{k|n} \tilde P_k$ with $\tilde P_1=\tilde Q_1$. 

\begin{conjecture}{Galois Groups for Centers of Hyperbolic Components}
\label{Conj:GaloisGroupsCenters} \lineclear
All $\tilde P_n$ are irreducible over $\Q$, and their Galois groups are the full symmetric groups.
\end{conjecture}
This would say that the centers of hyperbolic components of fixed period $n$ have the maximal symmetry possible. Manning confirmed this conjecture for low values of $n$ by computer experiments (unpublished).

Our approach for proving Theorem~\ref{Thm:GaloisGroups} will be using analytic continuation, like in the proof of the Ruffini-Abel theorem. This will yield explicit paths along which analytic continuation yields any given permutations. For specific values of $c\in\C$, the $P_n$ are polynomials in $\C[z]$ and factor over $\C$; we write them as $P_n(c)$.  Let 
\[
X_n:=\{c\in\C\colon \mbox{ all roots of $P_n(c)$ are simple}\} \,\,.
\]
Then all periodic points of period $n$ can be continued analytically through $X_n$, so the fundamental group of $X_n$ (with respect to any basepoint) acts on periodic points by analytic continuation. The question is which permutations can be so achieved. Of course, all permutations have to commute with the dynamics: if $z$ is a periodic points of $p_c$, then any permutation $\pi$ that is achieved by analytic continuation must have the property that $p_c(\pi(z))=\pi(p_c(z))$. It turns out that this is the only condition.


\begin{theo}{Analytic Continuation of Periodic Points}
\label{Thm:PermutationPeriodicPoints} \lineclear
For every period $n\ge 1$, analytic continuation in $X_n$ induces all permutations among periodic points of exact period $n$ that commute with the dynamics.
\end{theo}

\looseness -1
If $z_0$ is a double root of $P_n(c_0)$ for some $c_0\in\C$ , then $z_0$ is also a double root of $Q_n$ and $(d/dz)Q_n(z_0)=0$, hence $\mu:=(d/dz) p_{c_0}^{\circ n}(z_0)=1$. Here $\mu$ is the multiplier of the periodic orbit containing $z_0$. It is well known that a quadratic polynomial can have at most one non-repelling cycle. If $\mu=1$, then the indifferent orbit is called \emph{parabolic}; it is possible that the exact period of this orbit divides $n$, and the first return map of the parabolic orbit has a multiplier that may be a root of unity. Then $c_0$ is the root of a hyperbolic component of period $n$ (or dividing $n$) and thus the landing point of two parameter rays $R(\theta_1)$ and $R(\theta_2)$ with angles $\theta_i=a_1/(2^n-1)$. 
It follows that $\C\sm X_n$ is finite, and all periodic points of period $n$ can be continued analytically along curves in $X_n$ (as roots of $P_n(c)$). (However, $\bigcup_n \ovl{\C\sm X_n}=\partial\M$: every  $c\in\partial\M$ is a limit point of centers of hyperbolic components \cite{Orsay}, and it follows easily that the same is true for parabolics because for every $\eps>0$ almost every center has a parabolic parameter at distance less than $\eps$). 

Hyperbolic components come in two kinds, \emph{primitive} and \emph{satellite}, depending on the local properties of their roots. The two cases are as follows.

\begin{lemma}{Parabolic Parameters and Their Local Dynamics}
\label{Lem:LocalParabolics} \lineclear
Suppose that $p_{c_0}$ has a parabolic orbit of exact period $n$ with multiplier $\mu_0$. Then $c_0$ is the root of a unique hyperbolic component, and there are two possibilities:
\begin{description}
\item[The Primitive Case]
$\mu_0=1$ and the parabolic orbit is the merger of two orbits of period $n$. A small loop in parameter space around $c_0$ interchanges these two orbits and leaves all other orbits of period $n$ invariant. The parameter $c_0$ is the root of a single hyperbolic component of period $n$, and not on the boundary of any other hyperbolic component.

\item[The Satellite Case]
$\mu_0\neq 1$ is a $q$-th root of unity for some $q\ge 2$ and the parabolic orbit is the merger of one orbit of period $n$ and another orbit of period $qn$. A small loop in parameter space around $c_0$ induces a cyclic permutation on the period $qn$ orbit and leaves all other periodic orbits invariant, including the period $n$ orbit that becomes parabolic at $c_0$.
This cyclic permutation has cycles of length $q$, and thus it operates transitively on the orbit if and only if $n=1$. The parameter $c_0$ is the root of a unique hyperbolic component of period $qn$ and on the boundary of exactly one further hyperbolic component, and this component has period $n$. 
\end{description} 
\end{lemma}
\proof
The first return map of a parabolic orbit has the local form $f(z)=z\mapsto \mu_0z+z^{k+1} +O(z^{k+2})$ with $\mu_0$ a root of unity and $k\ge 1$. 

We first discuss the case $\mu_0=1$. Then there are $k$ parabolic petals at each of the parabolic periodic points, and each of them has to attract a critical orbit; since in our case there is only a single critical orbit, we have $k=1$. So $f(z)-z$ has a double zero at $z=0$ and the parabolic orbit splits up under perturbations of $c_0$ into two orbits of period $n$. Let $U$ be a neighborhood of $c_0$ in parameter space, small enough so that there are no other parabolic parameters of period $n$ in $U$. 

Analytic continuation of all orbits of period $n$ is possible locally for all parameters in $U\sm\{c_0\}$. Two orbits become parabolic at $c_0$, the others can be continued globally throughout $U$. Hence a loop in $U\sm\{c_0\}$ around $c_0$ either interchanges the two near-parabolic orbits, or it keeps all period $n$ orbits invariant. But in the latter case, the Open Mapping Principle implies that both of these orbits can become attracting near $c_0$, so $U$ intersects two hyperbolic components (or a single hyperbolic component that has $c_0$ twice on its boundary). But no two hyperbolic components of equal period touch \version{Proposition~\ref{Prop:Bifurcations}}{\cite{Orsay}, \cite[Sec.~6]{MiOrbits} or \cite[Corollary~5.7]{ExtRayMandel}}, and the boundary of any hyperbolic component is a Jordan curve \cite{Orsay}, \cite[Corollary~5.4]{ExtRayMandel}. Hence the two orbits are indeed interchanged. 

Since all orbits of periods other than $n$ are repelling for the parameter $c_0$, the point $c_0$ is not on the boundary of any hyperbolic component of period other than $n$, but of course it is on the boundary of a single hyperbolic component of period $n$.

If $\mu_0=e^{2\pi ip/q}$, then $p_{c_0}$ has a single parabolic orbit of period $n$, and it can be continued analytically in a neighborhood of $c_0$. The $q$-th iterate has the form $f^{\circ q}(z)=z+z^{k'+1}+O(z^{k'+2})$ with $k'$ petals at each parabolic point. Since $f$ relates groups of $q$ of these petals on one orbit, we have $k'=qm$, and since each orbit of petals has to attract a critical orbit as above, we have $m=1$ and $k'=q$. Under perturbation, the parabolic fixed point of $f^{\circ q}$ splits up into $q+1$ fixed points: these are the fixed point of $f$ and a single orbit of period $q$. All other orbits are repelling.

A small loop around $c_0$ can thus only act on the perturbed orbit of period $qn$, and it must do so by some cyclic permutation. The parameter $c_0$ is clearly on the boundary of at least one hyperbolic component of period $n$ and $qn$ each --- and since components of equal period never touch, it is on the boundary of exactly one component of period $n$ and $qn$, and not on the boundary of any other hyperbolic component.

Perturbing $c_0$ to nearby parameters $c$, every parabolic periodic point of $p_{c_0}$ breaks up into one point $w$ of period $n$ and $q$ points of period $qn$. These $q$ points form, to leading order, a regular $q$-gon with center $w$ because the first return map of $w$ has the local form $z\mapsto \mu z$ with $\mu$ near $e^{2\pi i p/q}$. When $c$ turns once around $c_0$, analytic continuation induces a cyclic permutation among these $q$ points of period $n$ so that during the loop, the $q$ points of period $n$ continue to lie on (almost) regular $q$-gons. When the loop is completed, the vertices of the $q$-gon are restored, so the $q$-gon will rotate by $s/q$ of a full turn, for some $s\ge 1$. If $s>1$, then the period $qn$ orbit and its multiplier would be restored to leading order after $c$ has completed $1/s$-th of a turn, and since boundaries of hyperbolic components are smooth curves this would imply that $c_0$ was on the boundary of $s$ hyperbolic components of period $qn$, and this is impossible as above.
Therefore $s=1$.
\qed

The fundamental group of $X_n$ (with respect to any basepoint) acts on the set of periodic points of $p_c$ of period $n$ by analytic continuation. 
Set $X:=\C\sm(\M\cup\R^+)$: this is a simply connected subset of all $X_n$ and will be used as a ``fat basepoint'' for the fundamental group of $X_n$.

\begin{figure}[htbp]
\includegraphics[width=100mm,trim=0 200 0 200,clip]{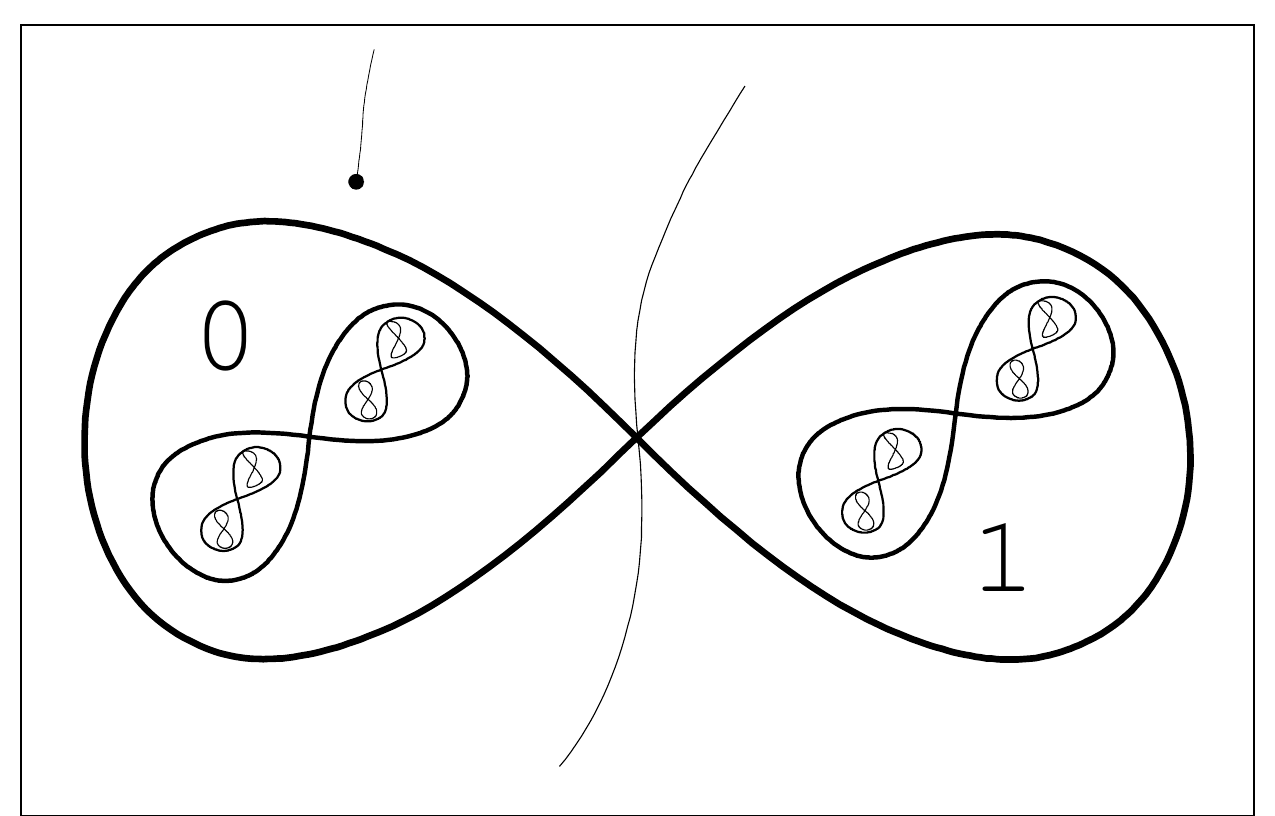}\\
\begin{picture}(00,0)
\put(-60,150){$c$}
\put(-59,170){$R_c(\theta)$}
\put(13,150){$R_c(\theta/2)$}
\put(-7,35){$R_c((\theta+1)/2)$}
\end{picture}
\caption{The partition for disconnected quadratic Julia sets is defined by the two dynamic rays that crash into the critical point; itineraries are defined with respect to this partition. The figure-8-curves are equipotentials (curves of equal speed of escape).}
\label{Fig:Itinerary}
\end{figure}
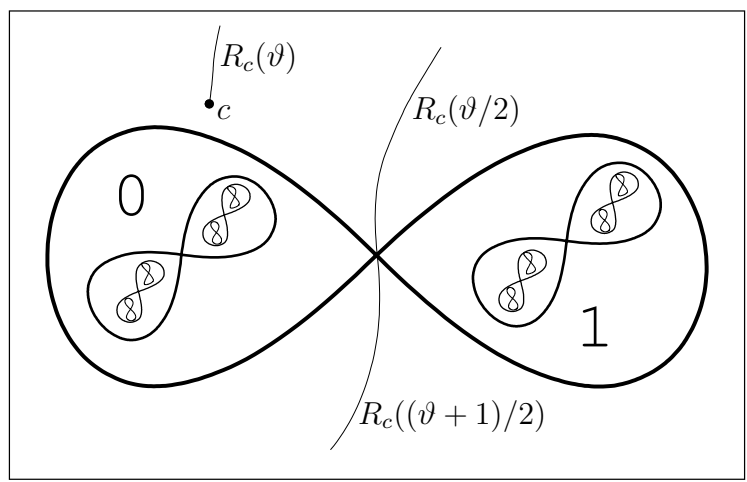

For every $c\in X$ we will describe periodic points of $p_c$ using symbolic dynamics: since $c\not\in\M$, the critical value $c$ is on the dynamic ray $R_c(\theta)$ for some $\theta\in\Circle$. Therefore the two dynamic rays $R_c(\theta/2)$ and $R_c((\theta+1)/2)$ both land at $0$ and separate the complex plane into two open parts, say $U_\0$ and $U_\1$ so that $c\in U_\1$ (see Figure~\ref{Fig:Itinerary}). The partition boundary does not intersect the Julia set $J_c$ of $p_c$, so we have $J_c\subset U_0\cup U_1$. Every $z\in J_c$ has an associated \emph{itinerary} $\tau_1\tau_2\tau_3\dots$, where $\tau_k\in\{\0,\1\}$ so that $p_c^{\circ {(k-1)}}(z)\in U_{\tau_k}$.

\begin{lemma}{Permutations and Symbolic Dynamics}
\label{Lem:PermutationsSymDyn} \lineclear
Let $c_0$ be the landing point of the parameter ray $R(\theta)$, where $\theta=a/(2^n-1)$. We consider the action of analytic continuation along a small loop around $c_0$ starting and ending at $R(\theta)$. 
\begin{itemize}
\item
If $c_0$ is the root of a primitive component, then analytic continuation along this loop interchanges the two periodic points with itineraries $\A(\nu(\theta))$ and $\Abar(\nu(\theta))$.
\item
If $c_0$ is the root of a satellite component, then analytic continuation along this loop affects the orbit of the point with itinerary $\Abar(\nu(\theta))$.
\end{itemize}
\end{lemma}
\proof
Let $R(\theta')$ be the second parameter ray landing at $c_0$ and let $U\subset\C$ be a disk neighborhood of $c_0$ in which $c_0$ is the only puncture of $X_n$, and which intersects no parameter ray at $n$-periodic angles other than $R(\theta)$ and $R(\theta')$. Let $U':=U\sm (R(\theta)\cup R(\theta'))$; this set is connected because $c_0\in U'$. 

Since the landing point of the periodic dynamic ray $R_c(\theta)$ depends continuously on the parameter whenever the ray lands (\cite[Expos\'e~XVII]{Orsay}, \cite[Proposition~5.2]{ExtRayMandel}), which is everywhere in $\C$ except on parameter rays at angles $R(2^k\theta)$ with $k\in\Nzero$, 
we can define a continuous function $z(c)$ as the landing point of the dynamic ray $R_c(\theta)$ for $c\in U'$. 

In the primitive case, analytic continuation along a simple loop around $c_0$ interchanges $z(c)$ with some point $z'(c)$ on a different orbit (Lemma~\ref{Lem:LocalParabolics}, using the fact that $R_{c_0}(\theta)$ lands on the parabolic orbit). The itinerary of $z$ can be determined most easily when $c\in U'$ is on a parameter ray $R(\tilde\theta)$ with $\tilde\theta$ near $\theta$: the itinerary of $z$  equals the itinerary of the angle $\theta$ under angle doubling with respect to the partition $\Circle\sm\{\tilde\theta/2,(\tilde\theta+1)/2\}$. This itinerary has period $n$, and for $\tilde\theta$ sufficiently close to $\theta$ it equals $\lim_{\phi\nearrow\theta}\nu(\phi)$ or $\lim_{\phi\searrow\theta}\nu(\phi)$ (depending on which side of $\theta$ the angle $\tilde\theta$ is), and these are $\A(\nu(\theta))$ or $\Abar(\nu(\theta))$. 

In the satellite case, $z(c)$ is on the orbit of period $n$ whenever $c$ is outside of the wake bounded by the two parameter rays $R(\theta)$ and $R(\theta')$ (inside the wake, or for $c=c_0$ on the wake boundary, the dynamic ray $R_c(\theta)$ lands on a point of lower period). This is the orbit that is affected by analytic continuation along loops around $c_0$ (Lemma~\ref{Lem:LocalParabolics} again), and the itinerary equals as above the itinerary of $\theta$ under angle doubling with respect to the partition $\Circle\sm\{\tilde\theta/2,(\tilde\theta+1)/2\}$ for $\tilde\theta$ near $\theta$ but outside our wake; using the convention that $\theta<\theta'$, then the itinerary equals $\lim_{\phi\searrow\theta}\nu(\phi)=\A(\nu(\theta))$. (The other limit $\lim_{\phi\nearrow\theta}\nu(\phi)$ equals $\Abar(\nu(\theta))$, and its period equals the period of the orbit at which the ray $R_c(\theta)$ lands within the wake; this is a proper divisor of $n$.)
\qed

\begin{prop}{Symmetric Permutation Group on Orbits}
\label{Prop:SymmetricGroup} \lineclear
Analytic continuation in $X_n$ induces the full symmetric group on the set of periodic orbits of period $n$.
\end{prop}
\proof
The domain $X_n$ is the complement of the finite set of roots of components of period $n$, and a small loop around any of these roots affects at most two periodic orbits of period $n$; if it does affect two orbits, then both orbits are interchanged. The permutation group among the periodic orbits of period $n$ is thus generated by pair exchanges. As soon as it acts transitively, it is automatically the full symmetric group.

It thus suffices to show that any orbit of period $n$ can be moved to the unique orbit containing the itinerary $\ovl{\1\1\dots\1\1\0}$. In fact, it suffices to show the following: suppose a periodic point has an itinerary containing at least two entries $\0$ during its period; then it can be moved to a periodic point whose itinerary has one entry $\0$ fewer per period. Repeated application will bring any periodic point onto the unique orbit with a single $\0$ per period, i.e., onto the orbit containing the itinerary  $\ovl{\1\1\dots\1\1\0}$.

Now consider a periodic point $z$ of period $n$ and assume that its itinerary $\tau_z$ contains at least two entries $\0$ per period. Let $\tau$ be the maximal shift of $\tau_z$ (with respect to the lexicographic order), and let $\tau'$ be the same sequence in which the $n$-th entry (which is necessarily a $\0$ in $\tau$) is replaced by a $\1$, again repeating the first $n$ entries periodically. Then by Lemma~\ref{Lem:MaxShift} there is a narrow hyperbolic component $W$ with associated kneading sequence $\nu(W)=\tau$. Let $R(\theta)$ be a parameter ray landing at the root of $W$; then $\A(\nu(\theta))=\nu(W)$, so the $\*$-periodic sequence $\nu(\theta)$ coincides with $\tau$ and $\tau'$ for $n-1$ entries. The component $W$ is primitive: by the remark after Definition~\ref{Def:Narrow}, a narrow component that is not primitive must bifurcate from the period $1$ component, and it would then have internal address $1\IntAddr n$ and kneading sequence with a single entry $\0$ in the period. Let $c_0$ be the root of $W$. By Lemma~\ref{Lem:PermutationsSymDyn}, a small loop around $c_0$ interchanges the periodic points with itineraries $\tau$ and $\tau'$. This is exactly the statement we need: we found a loop along which analytic continuation turns $z$ into a periodic point whose itinerary has one entry $\0$ fewer per period.
\qed

\proofof{Theorem~\ref{Thm:PermutationPeriodicPoints}}
Analytic continuation can achieve only permutations that commute with the dynamics. To see that all of them can actually be achieved, it suffices to specify one loop in $X_n$ that permutes one orbit transitively and leaves all other orbits (of the same period) unchanged: together with transitive permutations of all orbits, this generates all permutations that commute with the dynamics.

Let $c_0$ be the landing point of the parameter ray $R(1/(2^n-1))$; it is the bifurcation point from the period $1$ component to a component of period $n$. According to Lemma~\ref{Lem:LocalParabolics}, a small loop around $c_0$ induces a transitive permutation on a single orbit of period $n$ and leaves all other orbits unchanged. This proves the claim.
\qed

\begin{coro}{Riemann Surface of Periodic Points}
\label{Cor:RiemannSurface} \lineclear
For every $n\ge 1$, the analytic curve 
\[
\{ (c,z)\colon c\in X_n \mbox{ and $z$ is a periodic under $p_c$ of exact period $n$} \}
\]
 is connected, i.e., it is a Riemann surface.
\end{coro}

These results can be extended to preperiodic points as follows~\cite{Heidrun}.

\begin{coro}{Permutations of Preperiodic Points}
\label{Cor:PermutationsPreperiodic} \lineclear
Consider the set of preperiodic points that take exactly $k$ iterations to become periodic of period $n$. 
For each fixed pair of positive integers $k$ and $n$, analytic continuation along appropriate curves in $\C$ achieves all permutations that commute with the dynamics.
\end{coro}
\proof
The parameter ray $R(\theta_{k,n})$ with $\theta_{k,n}=1/(2^k2^{(n-1)})$ lands at a Misiurewicz-Thurston-parameter $c_{k,n}$ for which the dynamic ray $R_{c_{k,n}}(\theta)$ lands at the critical value \cite{Orsay,ExtRayMandel}. We have $\nu_{k,n}:=\nu(\theta_{k,n})=\1\1\dots\1\,\ovl{\1\1\dots\1\0}=\1^k\,\ovl{\1^{n-1}\0}$. A small loop around $c_{k,n}$ interchanges the two preperiodic points with itineraries $\0\nu_{k,n}$ and $\1\nu_{k,n}$ (these are the two preimage itineraries of the critical values): every periodic point can be continued analytically in a sufficiently small neighborhood of $c_{k,n}$, and the same is true for all points on the backwards orbits of these points as long as taking preimages does not involve the critical value. It follows that small loops around $c_{k,n}$ interchange the preperiodic points with itineraries $\tau\0\nu_{k,n}$ and $\tau\1\nu_{k,n}$ for every finite sequence $\tau$ over $\{\0,\1\}$.

Consider a preperiodic point $z$ with itinerary $\tau\1\,\ovl{\1^{n-1}\0}$, where $\tau$ is an arbitrary string over $\{\0,\1\}$ of length $k-1$ (the entry after $\tau$ must be $\1$, or the periodic part in the itinerary would start earlier). If $\tau$ has at least one entry $\0$, there is a value $k'$ so that a small loop around $c_{k',n}$ turns the last entry $\0$ within $\tau$ into an entry $\1$. Repeating this a finite number of times, $z$ can be continued analytically into the preperiodic point with itinerary $\nu_{k,n}$.
Analytic continuation thus acts transitively on the set of preperiodic points with itineraries $\tau\1\,\ovl{\1^{n-1}\0}$, for all $2^{k-1}$ sequences $\tau$ of length $k-1$. Since this is achieved by pair exchanges, the full symmetric group on these points is realized.

Two preperiodic points $z,z'$ of $p_c$ are on the same grand orbit if $p_c^{\circ n}(z)=p_c^{\circ n'}(z')$ for some positive integers $n,n'$. In terms of symbolic dynamics, this is the case if they have the same period, and the periodic parts of their itineraries are cyclic permutations of each other. A permutation of preperiodic points of preperiod $k$ and period $n$ on the same grand orbit commutes with the dynamics if and only if it induces the same cyclic permutations on the periodic parts of the orbit. 

For the grand orbit containing the point with itinerary $\ovl{\1^{n-1}\0}$, all permutations that commute with the dynamics can thus be achieved by analytic continuation around Misiurewicz-Thurston parameters $c_{k,n}$ and the root of the hyperbolic component $1_{1/n}\IntAddr n$ (a loop around the latter induces a transitive cyclic permutation of the periodic orbit containing the periodic point with itinerary $\ovl{\1^{n-1}\0}$). 

Since analytic continuation induces the full symmetric group on the set of grand orbits, the claim follows.
\qed

Note that any permutation of preperiodic points of preperiod $k$ and period $n$ induces a permutation of preperiodic points of preperiod in $\{k-1,k-2,\dots,1,0\}$ and period $n$. Analytic continuation takes place in $X_n$ from which finitely many Misiurewicz-Thurston points are removed.


\heading{Remark on Higher Degree Unicritical Polynomials}
Analogous results can also be obtained for the families of \emph{unicritical polynomials}, parametrized in the form $z\mapsto z^d+c$ for $d\ge 2$. All our results have generalizations to these families, and analytic continuation makes it possible to achieve all permutations of periodic points that commute with the dynamics; for details, see \cite[Section~12]{IntAdr}. Note that it is a much stronger statement to say that all permutations can be achieved by analytic continuation in the one-dimensional space of unicritical polynomials, rather than in the full $d-1$-dimensional space of general degree $d$ polynomials.

For preperiodic unicritical polynomials of degree $d>2$, there is one more invariant that is preserved under analytic permutation: preperiodic itineraries have the form $\tau=\tau_1\tau_2\dots\tau_k\,\ovl{\tau_{k+1}\dots\tau_{k+n}}$ with  $\tau_i\in\{\0,\1,\dots,\dit-\1\}$. Analytic continuation can move the preperiodic point with this itinerary to any other preperiodic point with itinerary $\tau'=\tau'_1\tau'_2\dots\tau'_k\,\ovl{\tau'_{k+1}\dots\tau'_{k+n}}$ provided $\tau'_{k+n}-\tau'_k=\tau_{k+n}-\tau_k$ (modulo $d$): not only must the length of preperiod and period be preserved, but also the ``cyclic difference'' between the last preperiodic point and its image point one full period later: both points have the same image, and the cyclic order among all points with that image is preserved. All permutations can be achieved that are compatible with the dynamics and that respect this condition.

A related study was done by Blanchard, Devaney, Keen \cite{BDK}: analytic continuation in the shift locus of degree $d$ polynomials realizes all automorphisms of the shift over $d$ symbols (in the special case of $d=2$, this corresponds to a loop around $\M$, and this interchanges all entries $\0$ and $\1$ in itineraries; indeed, this is the only non-trivial automorphism of the $2$-shift).

\looseness -1
A simple space where not all permutations can be achieved is the space of quadratic polynomials, parametrized as $z\mapsto \lambda z (1-z)$ with $\lambda\in\C$: the two fixed points are $z=0$ and $z=1-1/\lambda$ and they cannot be permuted by analytic continuation. This is related to the fact that the $\lambda$-space is not a true parameter space; every affine conjugacy class of quadratic polynomials is represented twice: the $\lambda$-space is the double cover over the true parameter space (written as $z\mapsto z^2+c$) that distinguishes the two fixed points. Another example is the space $f_c\colon z\mapsto (z^2+c)^2+c$ of second iterates of quadratic polynomials. Fixed points of such maps may have period $1$ or $2$ for $z\mapsto z^2+c$; this yields obstructions for permutations of fixed points of $f_c$. 

\hide{
We know of no parameter space of rational maps that is a true parameter space in which analytic continuation cannot achieve all permutations that commute with the dynamics.
}

\version{}{
\medskip{\sc Acknowledgements.}
Much of this work goes back to joint discussions with Eike Lau and the joint preprint \cite{IntAdr}. We would like to express our gratitude to him. We also gratefully acknowledge helpful and interesting discussions over the years with Christoph Bandt, Henk Bruin, Dima Dudko, John Hubbard, Karsten Keller, Misha Lyubich,  John Milnor, Chris Penrose, and many others, and we thank Simon Schmitt for his help with the pictures.}

\version{\reminder{Include in Sec 5 a quotable theorem that kneading sequences are admissible iff they satisfy the condition!}}{}

%% file: RefsIntAddr.tex

%% file: IntAddr2012.bbl
\begin{thebibliography}{MM9}




\bibitem[BDK]{BDK}
Paul Blanchard, Robert Devaney, Linda Keen,
\emph{The dynamics of complex polynomials and automorphisms of the shift}, Invent. Math. {\bf 104} (1991), 545--180.

\bibitem[Bo]{Bousch}
Thierry Bousch, \emph{Sur quelques probl\`emes de la dynamique holomorphe},
Thesis, Universit\'e de Paris-Sud (1992).

\bibitem[BF]{BF}
Bodil Branner, Nuria Fagella, \emph{Homeomorphisms between limbs of the Mandelbrot set}.  J. Geom. Anal.  {\bf 9} 3  (1999), 327--390.

\bibitem[BKS1]{TreesExist} Henk Bruin, Alexandra Kaffl, Dierk Schleicher,
\emph{Existence of quadratic Hubbard trees}. Fundamenta mathematicae, \textbf{202} (2009), 251--279.

\bibitem[BKS2]{SymDyn} Henk Bruin, Alexandra Kaffl, Dierk Schleicher,
\emph{Symbolic Dynamics of Quadratic Polynomials}, monograph in preparation.

\bibitem[BS1]{AdmissKneading} Henk Bruin, Dierk Schleicher, 
\emph{Admissibility of kneading sequences and structure of Hubbard trees for quadratic polynomials}. {Journal of the London Mathematical Society} {\bf 78} 2  (2008), 502--522; 
 
\bibitem[BS2]{BS-Admiss-PosMeasure} Henk Bruin, Dierk Schleicher,
\emph{Bernoulli measure of complex admissible kneading sequences}. Manuscript, submitted. arXiv:1205.1756.

\bibitem[BT]{BuffTanLei} Xavier Buff, Tan Lei, 
\emph{The quadratic dynatomic curves are smooth and irreducible}. Manuscript (2011).

\bibitem[D1]{DoAngles} Adrien~Douady,
{\em Algorithms for computing angles in the Mandelbrot set,}
in: {Chaotic dynamics and fractals}. Acad. Press (1986).

\bibitem[D2]{DoCompacts} Adrien~Douady,
{\em Descriptions of compact sets in $\C$}, in:
Topological Methods in Modern Mathematics, Publish or Perish
(1993) 429--465.


\bibitem[DH]{Orsay} Adrien~Douady, John~Hubbard,
{\em \'Etudes dynamique des polyn\^omes comples
I \& II},
Publ. Math. Orsay. (1984-85) {\em(The Orsay notes)}.

\bibitem[Du]{DimaDecoThm} Dzmitry Dudko, \emph{The Decoration Theorem for Mandelbrot and Multibrot Sets}. Manuscript, submitted. arXiv:1004.0633 (2010).
    

\bibitem[DS]{DimaHomeos} Dzmitry Dudko, Dierk Schleicher,
\emph{Homeomorphisms of limbs of the Mandelbrot set}. {Proceedings of the American Mathematical Society} \textbf{140} 6 (2012), 1947--1956. 

\bibitem[Ha]{HaissinskyTuning} Peter Haissinsky, 
\emph {Modulation dans l'ensemble de Mandelbrot}.
Tan Lei (ed.), {\em The Mandelbrot set, theme and variations},
Cambridge University Press {\bf 274} (2000) 37--65.


\bibitem[HY]{HY}
John~Hubbard,
{\em Local connectivity of Julia sets and bifurcation loci:
Three theorems of J.-C. Yoccoz,}
In: Topological Methods in Modern Mathematics, Publish or Perish
(1993) 467--511.

\bibitem[Ka]{VirpiPaper} Virpi~Kauko, {\em Trees of visible
components in the Mandelbrot set}, Fund. Math. {\bf
164} (2000) 41--60.


\bibitem[LS]{IntAdr} Eike~Lau, Dierk~Schleicher,
{\em Internal addresses in the Mandelbrot set and irreducibility
of polynomials,}
Stony Brook Preprint {\bf 19} (1994).

\bibitem[Lv]{Lavaurs} Pierre Lavaurs,
{\em Une description combinatoire de l'involution d\'efinie par
$M$ sur les rationnels \`a d\'enominateur impair},
C. R. Acad. Sci. Paris, S\'erie I Math. {\bf 303} (1986),
143--146.


\bibitem[M1]{MiBook} John~Milnor,
{\em Dynamics in one complex variable. Introductory lectures,}
Third edition, Princeton University Press (2005).

\bibitem[M2]{MiOrbits} John~Milnor,
{\em Periodic orbits, external rays, and the Mandelbrot
set: an expository account}. Ast\'erisque {\bf 261} (2000)
277--333.

\bibitem[M3]{MiRenorm}
John Milnor, \emph{Local connectivity of Julia sets: expository lectures}.
In: Tan Lei (ed.), \emph{The Mandelbrot set, theme and variations},  
London Math. Soc. Lecture Note Ser. {\bf 274}, Cambridge Univ. Press, Cambridge (2000), 
 67--116.

\bibitem[MT]{MilnorThurston} John Milnor, William Thurston, \emph{On iterated maps of the interval}. 
Springer LNM \textbf{1342} (1988), 465--563.
     

\bibitem[MP]{MP}
Patrick Morton, Pratiksha~Patel, \emph{The Galois theory of periodic points of polynomial maps},
Proc. Lond. Math Soc. (3) {\bf 68} (1994), 225--263.

\bibitem[M\"u]{Heidrun} Heidrun M\"undlein,
\emph{Permutationen vorperiodischer Punkte bei Iteration unikritischer Polynome} (permutations of preperiodic points under iteration of unicritical polynomials). Diploma thesis, Technische Universit\"at M\"unchen (1997).


\bibitem[Pe]{Pe} Chris~Penrose,
{\em On quotients of shifts associated with dendrite Julia sets
of quadratic polynomials,}
Thesis, University of Coventry, (1994).

\bibitem[Ri]{Riedl}
Johannes Riedl, \emph{Arcs in Multibrot Sets, Locally Connected Julia Sets 
and Their Construction by Quasiconformal Surgery}.
PhD thesis, Technische Universit\"at M\"unchen (2000). 


\bibitem[RS]{CrossRenorm} Johannes Riedl, Dierk Schleicher, 
{\em On the locus of crossed renormalization}. In: {\em Proceedings of the Research Institute of Mathematical
Sciences}, Kyoto {\bf 4} (1998), 21--54.



\bibitem[Sch1]{ExtRayMandel} Dierk~Schleicher,
{\em Rational external rays of the Mandelbrot set},
Asterisque {\bf 261} (2000) 405--443.

\bibitem[Sch2]{Fibers2} Dierk~Schleicher,
{\em On fibers and local connectivity of Mandelbrot and Multibrot
sets}. In: M.~Lapidus, M.~van Frankenhuysen (eds): {\em
Fractal Geometry and Applications: A Jubilee of Benoit Mandelbrot}.
Proceedings of Symposia in Pure Mathematics {\bf 72}, 
American Mathematical Society (2004), 477--507.


\bibitem[T]{thurston} William Thurston,
{\em On the geometry and dynamics of iterated rational maps,}
in: \emph{Complex dynamics: families and friends}, ed. Dierk Schleicher, A K Peters, Wellesley, MA, 2009, 3--137.




\end{thebibliography}
